\documentclass[preprint,10pt]{elsarticle}

\biboptions{square,sort,comma,numbers}

\usepackage[latin1]{inputenc}
\usepackage{latexsym}
\usepackage{graphicx}

\usepackage{stmaryrd}

\usepackage{amscd}
\usepackage{amssymb,amsmath,amsthm}
\usepackage{mathptmx}
\usepackage{mathrsfs}
\usepackage{color}
\usepackage{xspace}
\usepackage{bussproofs}

\usepackage{tikz}
\usepackage{bpextra}

\usepackage{centernot}

\usepackage{array}
\newcolumntype{C}[1]{>{\centering\arraybackslash}p{#1}}

\usepackage{bigstrut}

\usepackage{hyperref}

\usepackage{cleveref}

\usepackage{relsize}

\usepackage{caption}

\allowdisplaybreaks[4]

\DeclareSymbolFont{letters}{OML}{cmm}{m}{it}

\DeclareMathAlphabet{\mathcal}{OMS}{cmsy}{m}{n}

\newtheorem{theorem}{Theorem}[section]

\newtheorem{lemma}[theorem]{Lemma}
\newtheorem{proposition}[theorem]{Proposition}

\newtheorem{claim}{Claim}
\newtheorem*{claim*}{Claim}

\theoremstyle{definition}
\newtheorem{definition}[theorem]{Definition}
\newtheorem{example}[theorem]{Example}

\newtheorem*{theorem*}{Theorem}

\newenvironment{proofclaim}[1][the claim]{\par
\noindent \emph{Proof of #1.} }
{\hfill$\dashv$\vspace{4pt}}

\usepackage{multirow}

\newcommand{\Inql}{\ensuremath{\mathbf{InqL}}\xspace}

\newcommand{\ML}{\ensuremath{\mathsf{ML}}\xspace}

\newcommand{\PD}{\ensuremath{\mathbf{PD}}\xspace}
\newcommand{\PID}{\ensuremath{\mathbf{PID}}\xspace}

\newcommand{\IPC}{\ensuremath{\mathbf{IPL}}\xspace}

\newcommand{\CPC}{\ensuremath{\mathbf{CPL}}\xspace}
\newcommand{\CPL}{\ensuremath{\mathbf{CPL}}\xspace}

\newcommand{\PDbor}{\ensuremath{\PD^{\bor}}\xspace}

\newcommand{\ND}{\ensuremath{\mathsf{ND}}\xspace}
\newcommand{\KP}{\ensuremath{\mathsf{KP}}\xspace}

\newcommand{\PT}{\ensuremath{\mathbf{PT}_0}\xspace}

\newcommand{\dep}{\ensuremath{\mathop{=\!}}\xspace}

\newcommand{\sor}{\ensuremath{\otimes}\xspace}
\newcommand{\bor}{\ensuremath{\vee}\xspace}

\DeclareMathOperator*{\bigsor}{\bigotimes}
\DeclareMathOperator*{\bigbor}{\bigvee}

\newcommand{\ci}{\ensuremath{\wedge\textsf{I}}\xspace}

\newcommand{\ce}{\ensuremath{\wedge\textsf{E}}\xspace}
\newcommand{\sori}{\ensuremath{\sor\textsf{I}}\xspace}

\newcommand{\sorwe}{\ensuremath{\sor\textsf{E}^-}\xspace}
\newcommand{\bote}{\ensuremath{\mathbf{\bot}\textsf{E}}\xspace}
\newcommand{\sors}{\ensuremath{\sor\textsf{Sub}}\xspace}
\newcommand{\bori}{\ensuremath{\bor\textsf{I}}\xspace}
\newcommand{\bore}{\ensuremath{\bor\textsf{E}}\xspace}
\newcommand{\exfalso}{\ensuremath{\it{\textsf{ex falso}}}\xspace}
\newcommand{\exclmid}{\ensuremath{\textsf{EM}_0}\xspace}

\newcommand{\com}{\ensuremath{\textsf{Com}}\xspace}
\newcommand{\ass}{\ensuremath{\textsf{Ass}}\xspace}
\newcommand{\dstr}{\ensuremath{\textsf{Dstr}}\xspace}

\newcommand{\cs}{\ensuremath{\wedge\textsf{Sub}}\xspace}

\newcommand{\boti}{\ensuremath{\bot\textsf{I}}\xspace}

\newcommand{\depiz}{\ensuremath{\textsf{DepI}_0}\xspace}
\newcommand{\depez}{\ensuremath{\textsf{DepE}_0}\xspace}

\newcommand{\depik}{\ensuremath{\textsf{DepI}_k}\xspace}
\newcommand{\depek}{\ensuremath{\textsf{DepE}_k}\xspace}
\newcommand{\se}{\ensuremath{\textsf{SE}}\xspace}

\newcommand{\dstrs}{\ensuremath{\textsf{Dstr}^\ast\!}\xspace}

\newcommand{\MP}{\ensuremath{\mathsf{MP}}\xspace}
\newcommand{\Sub}{\ensuremath{\mathsf{Sub}}\xspace}

\definecolor{darkgreen}{rgb}{0,.7,0}

\usepackage[normalem]{ulem}

\newcommand{\LL}{\ensuremath{\mathsf{L}}\xspace}

\journal{arXiv.org}

\begin{document}

\begin{frontmatter}



\title{Propositional Logics of Dependence\tnoteref{tn}}

\tnotetext[tn]{Most of the results in this paper were included in the dissertation of the first author \cite{Yang_dissertation}, which was supervised by the second author.}

 \author[FY]{Fan Yang\corref{cor}\fnref{fnf}}

  \address[FY]{Department of Philosophy and Religious Studies, Utrecht University, Janskerkhof 13, 3512 BL Utrecht, The Netherlands}
\ead{fan.yang.c@gmail.com}

\fntext[fn]{The research was carried out in the Graduate School of Mathematics and Statistics of the University of Helsinki.}

 \author[JV]{Jouko V\"{a}\"{a}n\"{a}nen\fnref{fnj}}

  \ead{jouko.vaananen@helsinki.fi} 
 \address[JV]{Department of Mathematics and Statistics, Gustaf H\"{a}llstr\"{o}min katu 2b, PL 68, FIN-00014 University of Helsinki, Finland and University of Amsterdam, The Netherlands}

\fntext[fnj]{The research was partially supported by grant 251557 of the Academy of Finland.}

\cortext[cor]{Corresponding author}



\begin{abstract}

In this paper, we study logics of dependence on the propositional level. We prove that several interesting propositional logics of dependence, including \emph{propositional dependence logic}, \emph{propositional intuitionistic dependence logic} as well as \emph{propositional inquisitive logic}, are expressively complete and have disjunctive or conjunctive normal forms. We provide deduction systems and prove the completeness theorems for these logics.

\end{abstract}


\end{frontmatter}


\section{Introduction}\label{sec:intro}


The idea of dependence logic, introduced in \cite{Van07dl} on the basis of Hodges \cite{Hodges1997b}, is the following: If truth in the traditional sense, i.e. as  Tarski defined it,  is defined with respect to a {\em set} of assignments, rather than just {\em one} assignment, it becomes possible to talk meaningfully about variables being {\em dependent} or {\em independent} from each other. The set of assignments can be thought of as a {\em data set}, giving evidence---much as in statistics---about mutual dependencies between the variables. If the set of assignments---the data---is thought of as a {\em database} we arrive at the concept of {\em functional dependency}, an important concept in database theory  since Codd's pioneering paper \cite{DBLP:persons/Codd71a}.
We can also consider the set of assignments as expressing {\em uncertainty} about one ``true" assignment, as in {\em inquisitive logic } (\cite{InquiLog}), or as indicating {\em belief} about an unknown assignment, as in {\em doxastic logic} (\cite{Galliani2013}). Finally, considering truth as given by sets of assignments leads naturally to the concept of the {\em probability}   that a randomly chosen assignment (from the given set) satisfies a given propositional formula. This idea is developed in \cite{HPV} to analyse so-called Bell's Inequalities of Quantum Physics. 
 
Following \cite{Van07dl}, we call sets of assignments {\em teams}. Teams have been previously used to study  dependence concepts in predicate logic \cite{Van07dl} and modal logic \cite{VaMDL08}. We now  focus this study to propositional logic. The fundamental concept of dependence logic is the concept $\dep(\vec{x},y)$ of a variable $y$ depending on a sequence $\vec{x}$ of other variables, which is taken as a new atomic formula. The meaning of such atomic dependence formulas is given via  {\em teams}. 

Studying the logics of dependence concepts in propositional logic resembles the case of predicate logic in that we use the method of teams. A \emph{team} in this case is defined to be \emph{a set of valuations} of propositional variables. There are, however, also significant differences between the predicate logic and the propositional cases. Notably, propositional logics of dependence are {\em decidable}. This is because for any given formula of the logics with $n$ propositional variables, there are in total $2^n$ valuations and $2^{2^n}$ teams.  The method of truth tables has its analogue in these logics, but the size of such tables grows exponentially faster than in the case of classical propositional logic, rendering it virtually inapplicable. This emphasizes the role of the axioms and the completeness theorem in providing a manageable alternative for establishing logical consequence.

Classical propositional logic is based on propositions of the form 
\begin{center}
\begin{tabular}{l}
$p$\\
Not $p$\\
$p$ or $q$\\
If $p$, then $q$
\end{tabular}
\end{center}
and more generally
\begin{equation}\label{cpc_imp_intuition}
\begin{array}{l}
\mbox{If $p_{i_1},\ldots,p_{i_k}$, then $q$}.\\
\end{array}
\end{equation}
We present extensions of classical propositional logic in which one can express, in addition to the above,  propositions of the form
``$q$ depends on $p$",
or more generally 
\begin{equation}\label{dep_intuition}
\mbox{$q$ depends on $p_{i_1},\ldots,p_{i_k}$.}
\end{equation}
In our setting, (\ref{dep_intuition})  is treated as an atomic fact. This is expressed formally by a new atomic formula
\begin{equation}\label{dep_atom_example5}
\dep(p_{i_1},\ldots,p_{i_k},q),
\end{equation}
which we call the \emph{dependence atom}.

 Intuitively, (\ref{dep_intuition}) means that to know whether $q$ holds it is sufficient to consult the truth values of  $p_{i_1},\ldots,p_{i_k}$. Note that, as in the first-order dependence logic case,  (\ref{dep_intuition})  says nothing about the {\em way} in which $p_{i_1},\ldots,p_{i_k}$ are logically related to $q$. It may be that $p_{i_1}\wedge\ldots\wedge p_{i_k}$ logical implies $q$, or that $\neg p_{i_1}\wedge\ldots\wedge \neg p_{i_k}$ logical implies $\neg q$, or anything in between. Technically speaking, this is to say:
\begin{equation}\label{dep_meaning}
\mbox{The truth value of $q$ is a function of the truth values of $p_{i_1},\ldots,p_{i_k}$.}
\end{equation}


Some examples of natural language sentences involving dependence are the following:
\begin{enumerate}
\item \it{Whether it rains depends completely on whether it is winter or summer.}
\item Whether you end up in the town depends entirely on whether you turn  left here or right.
\item I will be absent depending on whether he shows up or not.
\end{enumerate}

We now define propositional dependence logic formally. We start by recalling classical propositional logic and define our {\em team semantics} for it. Although team semantics is intended for the extension of classical propositional logic obtained by adding the new atomic formulas, namely the dependence atoms of (\ref{dep_atom_example5}), it has some interesting uses in the classical case, too, as we will see below.



\begin{definition}\label{syntax_cpc}
Let $\{p_i\mid i\in\mathbb{N}\}$ be the set of propositional variables.
Well-formed
formulas of \emph{classical propostional logic} (\CPC)  are given by the following grammar
\[
    \phi::= \,p_i\mid \neg p_i\mid\bot\mid(\phi\wedge\phi)\mid(\phi\sor\phi)\mid (\phi\to\phi).
\] 
\end{definition}
We use the symbol $\otimes$ to denote the disjunction of \CPC and reserve the usual disjunction symbol $\vee$ for another use. Hereafter we adopt the usual convention for negation and write $\neg\phi$ for the formula $\phi\to\bot$. We will see  from the definition of team semantics below that the negated propositional variable $\neg p_i$ has the same semantics as $p_i\to \bot$.

We now give the formal definition of the crucial notion of the team semantics, namely \emph{teams}.

\begin{definition}  A \emph{valuation} is a function from the set $\mathbb{N}$ of natural numbers to the set $2=\{0,1\}$. A \emph{team} $X$ is a set of valuations. That is, $X\subseteq 2^{\mathbb{N}}$.
\end{definition}
See \Cref{rer} for an example of a team.  One can think of a team such as that given in \Cref{rer} as the result of five {\em tests} about the atomic propositions $p_1,p_2,...$. The atomic propositions may code basic data about people and the five assignments $s_1,...,s_5$ may be this data for five different people. The atomic propositions $p_1,p_2,...$ may also record moves in a game. The five rows represent then five plays of the same game. There are many possible intuitions about teams and we do not fix any particular intuition. The empty set $\emptyset$ and the full set $2^{\mathbb{N}}$ are special cases of teams. Another rather special case is the {\em singleton team} $\{s\}$. The usual semantics of propositional logic resembles team semantics with singleton teams. It is the non-singleton teams that bring out the new phenomena of team semantics.

\begin{table}[t]
\begin{center}
\begin{tabular}{c|cccccc}
&$p_1$&$p_2$&$p_3$&$p_4$&$p_5$&$\dots$\\
\hline
$s_1$&1&1&1&0&1\\
$s_2$&1&0&1&0&0\\
$s_3$&1&1&1&1&1&$\dots$\\
$s_4$&1&0&1&1&0\\
$s_5$&1&0&0&1&1\\
\end{tabular}
\caption{A team $X=\{s_1,s_2,s_3,s_4,s_5\}$}\label{rer}
\end{center}
\end{table}%

\begin{definition}\label{TS_cpc}
We inductively define the notion of a formula $\phi$ in the language of \CPC  being \emph{true} on a team $X$, denoted $X\models\phi$, as follows:
\begin{itemize}
\item $X\models p_i$ iff
for all $s\in X$, $s(i)=1$
  \item $X\models\neg p_i$  iff
for all $s\in X$, $s(i)=0$
  \item $X\models\bot$ iff $X=\emptyset$

  \item $X\models\phi\wedge\psi$ iff $X\models\phi$ and
  $X\models\psi$
  \item $X\models\phi\sor\psi$ iff there exist two subteams $Y,Z\subseteq X$ with $X=Y\cup Z$ such that 
  \(Y\models\phi\text{ and }Z\models\psi\)
  \item $X\models \phi\to\psi$ iff for any subteam $Y\subseteq X$, 
  \(Y\models \phi\,\Longrightarrow \,Y\models\psi\)
\end{itemize}
\end{definition}


In this paper, we also call formulas in the language of \CPL \emph{classical formulas}. An easy inductive proof shows that classical formulas $\phi$ are \emph{flat}, that is,
\begin{description}
\item[(Flatness Property)] $X\models\phi\iff\forall s\in X(\{s\}\models\phi)$ holds for all teams $X$.
\end{description}
This means that the truth of a classical formula under the team semantics is determined by its truth over singleton teams.
Over singleton teams the disjunction $\otimes$ and the implication $\to$ behave classically:
$$\{s\}\models\phi\sor\psi\iff\{s\}\models\phi\mbox{ or }\{s\}\models\psi$$
$$\{s\}\models\phi\to\psi\iff\{s\}\not\models\phi\mbox{ or }\{s\}\models\psi.$$
From this one can easily show that for classical formulas the team semantics over singleton teams coincides with the usual truth semantics of \CPL, that is,
\begin{equation}\label{cpl_team_red}
\{s\}\models\phi\iff s\models\phi \text{ holds for all classical formulas }\phi
\end{equation}
where $s\models\phi$ is meant in the ordinary sense of a valuation $s$ satisfying $\phi$. In this sense, the team semantics for classical formulas actually reduces to the usual single-valuation semantics. 



Let $|\phi|_X=\{s\in X : \{s\}\models\phi\}$. For any classical formula $\phi$ and any team $X$
$$\begin{array}{lcl}
X\models\phi&\mbox{ implies }&|\phi|_X=X\\
X\models\neg\phi&\mbox{ implies }&|\phi|_X=\emptyset\\
\end{array}$$ but the set $|\phi|_X$ can also be a proper non-empty subset of $X$, e.g., the set $|p_2|_X=\{s_1,s_3\}$ of the team $X$ of \Cref{rer}. In \cite{HPV} the number $|\,|\phi|_X|/|X|$ (for finite non-empty $X$) is for obvious reasons called the {\em probability} of $\phi$ in $X$, denoted $[\phi]_X$. When such probabilities are considered it makes sense to allow the same valuation to occur in $X$ many times (i.e., to consider multisets). In \cite{HPV} a complete axiomatization is given for propositions about linear combinations of probabilities $[\phi]_X$. However, in this paper we proceed in a different direction and use teams to define truth values for dependence statements. We now define the logic for expressing dependence statements.


\begin{definition}\label{syntax_pt}
 Well-formed
formulas of \emph{propositional downward closed team logic} (\PT)  are given by the following grammar
\[
    \phi::= \,p_i\mid \neg p_i\mid\bot\mid\dep(p_{i_1},\dots,p_{i_k},p_j)\mid(\phi\wedge\phi)\mid(\phi\sor\phi)\mid (\phi\vee\phi)\mid (\phi\to\phi)
\] 
\end{definition}

We call the connectives $\sor$, $\vee$ and $\to$  \emph{tensor (disjunction)},  \emph{intuitionistic disjunction} and   \emph{intuitionistic implication}, respectively.

\begin{definition}\label{TS_PT}
We inductively define the notion of a formula $\phi$ in the language of \PT  being \emph{true} on a team $X$, denoted $X\models\phi$. All the cases are identical to those defined in \Cref{TS_cpc} and additionally:
\begin{itemize}
  \item $X\models\,\dep(p_{i_1},\dots,p_{i_k},p_j)$ iff for all $s,s'\in X$,
  \[[\,s(i_1)=s'(i_1),\,\dots,\,s(i_{k})=s'(i_{k})\,]~\Longrightarrow ~s(j)=s'(j);\]

  \item $X\models \phi\bor\psi$ iff $X\models \phi$ or $X\models\psi$;
\end{itemize}

We say that a formula $\phi$  in the language of \PT is \emph{valid}, denoted  $\models\phi$, if $X\models\phi$ holds for all teams $X$. We write $\phi\models\psi$ if $X\models\phi\Longrightarrow X\models\psi$ holds for all teams $X$, and say then that $\psi $ is a \emph{semantic consequence} of $\phi$.
 If $\phi\models\psi$ and $\psi\models\phi$, then we say that $\phi$ and $\psi$ are \emph{ semantically equivalent}, in symbols $\phi\equiv\psi$. 
\end{definition}

 
Immediately from the team semantics we obtain that formulas in the language of \PT have the properties stated in the following theorem. These properties  also hold in the case of  first-order dependence  logic (see e.g. \cite{Van07dl}).
\begin{theorem}\label{basic_prop_pt}
Let $\phi$ be a formula in the language of \PT, and $X$ and $Y$ two teams.  
\begin{description}
\item[(Locality Property)] If $\{p_{i_1},\dots,p_{i_n}\}$ is the set of propositional variables occurring in $\phi$ and 
$\{s\upharpoonright \{{i_1},\dots,{i_n}\}: s\in X\}=\{s\upharpoonright \{{i_1},\dots,{i_n}\}:s\in Y\}$, then
\[X\models\phi\iff Y\models\phi.\] 
\item[(Downward Closure Property)] If $X\models\phi$ and $Y\subseteq X$, then $Y\models\phi$.
\item[(Empty Team Property)] $\emptyset\models\phi$.
\end{description}
\end{theorem}


Let us spend a few words on the  atoms and the connectives of our logic \PT. To get an idea  of the definition of $X\models\,\dep(p_{i_1},\dots,p_{i_k},p_j)$ let us consider an example. In \Cref{rer} we have a team $X$ satisfying the dependence atom   $\dep(p_2,p_3,p_5)$. This means that if we take any two rows of the team with the same truth values for both $p_2$ and $p_3$, then also the truth value of $p_5$ is the same. On rows 1 and 3 we have $p_2$ and $p_3$ true, and indeed on both rows $p_5$ has the same truth value, namely true. Similarly, on rows 2 and 4 we have on both rows $p_2$ false and $p_3$ true, and indeed on both rows $p_5$ has the same truth value, namely false. This is often (especially in database theory) called {\em functional dependency}, because there is a function which gives the value of $p_5$ with the values of $p_2$ and $p_3$ as arguments. In this case the function is 
\[(1,1)\mapsto 1,\quad
(0,1)\mapsto 0,\quad
(0,0)\mapsto 1.\]
Functional dependencies are axiomatized in database theory by the so-called {\em Armstrong's Axioms \cite{Armstrong_Axioms} as follows\footnote{In this paper we only consider a simplified version of Armstrong's Axioms. In full generality  Armstrong's Axioms have arbitrary finite sequences of variables.} (presented in terms of dependence atoms)}:
\begin{description}
\item[(i)] $\dep(x,x)$
\item[(ii)] if $\dep(x,y,z)$, then $\dep(y,x,z)$
\item[(iii)] if $\dep(x,x,y)$, then $\dep(x,y)$
\item[(iv)] if $\dep(y,z)$, then  $\dep(x,y,z)$
\item[(v)] if $\dep(x,y)$ and $\dep(y,z)$, then  $\dep(x,z)$
\end{description}
But these rules do not characterize the interaction between dependence atoms and other formulas. We will define in \Cref{arm} natural deduction rules for dependence atoms that also characterize such interactions, and give a simple derivation of Armstrong's Axioms from those rules. 

One distinguishing feature of dependence atoms is that they are in general \emph{not} flat. To see why, consider an arbitrary dependence atom $\dep(p_{i_1},\dots,p_{i_k},p_j)$ and an arbitrary valuation $s$. By the semantics we know that $\{s\}\models\dep(p_{i_1},\dots,p_{i_k},p_j)$ is always true (for $\{s\}$ is a singleton). However, obviously not every dependence atom is true on every team. For example, the team $X$ from \Cref{rer} falsifies $\dep(p_3,p_4,p_5)$, as, e.g., $s_1(p_3)=1=s_2(p_3)$ and $s_1(p_4)=0=s_2(p_4)$, whereas $s_1(p_5)=1\neq0= s_2(p_5)$.

The negation of an arbitrary formula $\phi$ in the language of \PT is defined as $\neg\phi:=\phi\to\bot$. 
Such defined negations $\neg\phi$ are always flat. Since dependence atoms are true on singleton teams, the negation  of a dependence atom (a flat formula) is true only on the empty team, i.e., $X\models\neg\dep(p_{i_1},\dots,p_{i_k},p_j)$ iff $X=\emptyset$, or $\neg\dep(p_{i_1},\dots,p_{i_k},p_j)\equiv\bot$. For other interesting properties of the negation of team semantics we refer the reader to \cite{IemhoffYang15}.


The semantics of tensor $\sor$ resembles the semantics of tensor in linear logic; see \cite{AbVan09} for further discussion.  To understand the semantics of $\phi\sor\psi$, we may think of a team $X$ as representing uncertainty about an unknown valuation $s$, then the representation $X=Y\cup Z$ narrows the uncertainty into two more certain cases $Y$ and $Z$ that respectively satisfy each of the disjuncts. 
Since the tensor $\sor$ of \PT originates from the disjunction of \CPL, it is not surprising (and we leave it for the reader to check) that the \emph{Law of Excluded Middle} with respect to $\sor$ holds for classical formulas $\alpha$, i.e., $\models\alpha\sor\neg\alpha$, and, in fact, it holds also for all flat formulas. For non-flat formulas the Law fails in general; e.g., $\not\models \dep(p_i)\sor\neg \dep(p_i)$.

The set $\mathcal{P}(2^{\mathbb{N}})\setminus\{\emptyset\}$ of all nonempty teams endowed with the superset relation $\supseteq$ forms an intuitionistic Kripke frame $\mathfrak{F}=(\mathcal{P}(2^{\mathbb{N}})\setminus\{\emptyset\},\supseteq)$. The team semantics for the intuitionistic disjunction $\vee$ and the intuitionistic implication $\to$ resembles the Kripke semantics of the disjunction and the implication of intuitionistic propositional logic over the fixed Kripke frame $\mathfrak{F}$. We will come back to this issue in \Cref{sec:axiom_pid}.
For now we want to point out that the intuitionistic disjunction of \PT does have a constructive feature, as it satisfies the \emph{Disjunction Property}:
\begin{description}
\item[(Disjunction Property)] 
If  $\models\phi\bor\psi$, then $\models \phi$ or $\models \psi$. 
\end{description}
To see why this holds, 
let $2^{\mathbb{N}}$ be the biggest team.
By the downward closure property we have  
\[\models\phi\bor\psi\Longrightarrow~ 2^{\mathbb{N}}\models \phi\bor\psi\Longrightarrow~ 2^{\mathbb{N}}\models \phi\text{ or }2^{\mathbb{N}}\models \psi\Longrightarrow~ \models\phi\text{ or }\models\psi.\]
For the intuitionistic implication of \PT the \emph{Deduction Theorem} holds:
\begin{description}
\item[(Deduction Theorem)] 
\(\Gamma,\phi\models\psi\text{ if and only if }\Gamma\models\phi\to\psi.\)
\end{description}
Over singleton teams these two intuitionistic connectives behave in the classical manner. We already discussed this for the intuitionistic implication (which has the same team semantics as the implication of \CPL); for the intuitionistic disjunction we have 
\[\{s\}\models \phi\vee\psi\iff \{s\}\models \phi\text{ or }\{s\}\models \psi.\]



 \begin{table}[t]
\centering
\begin{tabular}{c|cccc|cccc}
&$p_{i_1}$&$p_{i_2}$&$\cdots$&$p_{i_n}$&$X_1$&$X_2$&$\cdots$&$X_{2^{2^n}}$\\\hline
$s_1$&$1$&$1$&$\cdots$&$1$&$1$&$1$&$\cdots$&1\\
$s_2$&$0$&$1$&$\cdots$&$1$&$0$&$1$&$\cdots$&1\\
$\vdots$&&&$\vdots$&&$\vdots$&&$\vdots$&\\
$s_{2^n}$&$0$&$0$&$\cdots$&$0$&$0$&$0$&$\cdots$&1\\
\hline
$\phi$&&&&&$1$&$0$&$\cdots$&$1$\\
\end{tabular}
\caption{A truth table for a formula $\phi(p_{i_1},\dots,p_{i_n})$}
\label{truth_table}
\end{table}%

In classical propositional logic \CPL a fundamental and intuitive method to study the meaning of propositional formulas is the {\em truth table method}. The truth table methodically lists all the possible valuations and a simple computation yields $1$ or $0$ as the truth value of the formula for {\em that} valuation. Once the table is drawn, a column of $1$'s for $\phi$ means the formula $\phi$ is valid and a column with at least one $1$ means $\phi$ is satisfiable. For a formula of $n$ propositional variables the table has $2^n$ rows, whence this method is said to be of {\em exponential} (time) complexity, but conceived as a guessing problem the satisfiability problem is \mbox{NP}-complete. Respectively, the validity problem is \mbox{co-NP}-complete. When we move to propositional downward closed team logic \PT, the situation is one step more complicated. We still have the {\em truth table method} but now the truth table has two parts. In \Cref{truth_table} we have on the left a table of all possible valuations and the last columns ``$X$'' make a choice which valuations are taken to the team $X$. Finally, if a formula $\phi(p_{i_1},\dots,p_{i_n})$ is given, a systematic list of its truth values in different teams can be listed at the bottom of the table. All this is in principle doable, but of course the size of the table quickly explodes.

Let $N$ be an $n$-element set of natural numbers. We call a function $s:N\to \{0,1\}$ a \emph{valuation on} $N$, and a set of valuations on $N$ a \emph{team on} $N$. For example, the teams $X_i$ in \Cref{truth_table} are teams on $\{i_1,\dots,i_n\}$. There are in total $2^n$ distinct valuations on $N$, and $2^{2^n}$ distinct teams on $N$, among which there is a biggest team (denoted  $2^N$) consisting of all valuations on $N$. For any formula $\phi$ in the language of \PT we write $\phi(p_{i_1},\dots,p_{i_n})$ to mean that the propositional variables occurring in $\phi$ are among $p_{i_1},\dots,p_{i_n}$. By the locality property, to evaluate a formula $\phi(p_{i_1},\dots,p_{i_n})$ it is sufficient to consider teams on $N=\{i_1,\dots,i_n\}$ only. We write $\llbracket \phi\rrbracket$ for the set of all teams on $N$ satisfying $\phi$, i.e.,
\[\llbracket \phi\rrbracket:=\{X\subseteq 2^N\mid X\models\phi\}.\]
It follows from \Cref{basic_prop_pt} that the set $\llbracket \phi\rrbracket$ is nonempty (since $\emptyset\in\llbracket\phi\rrbracket$) and closed under subsets. Denote by $\nabla_N$ the family of all nonempty downward closed collections of teams on $N$, i.e., 
\[\nabla_N=\{\mathcal{K}\subseteq 2^{2^N}\mid \mathcal{K}\neq \emptyset\text{ and }({Y}\subseteq{X}\in \mathcal{K}\Longrightarrow {Y}\in\mathcal{K})\}.\]
We will prove in this paper that the logic \PT is \emph{expressively complete} in the sense that given any set $N=\{i_1,\dots,i_n\}$ of natural numbers, 
\begin{equation}\label{expres_compl_df}
\nabla_N=\{\llbracket \phi\rrbracket:~\phi(p_{i_1},\dots,p_{i_n})\text{ is a formula in the language of \PT}\}.
\end{equation}

We study in this paper also some fragments of \PT whose well-formed formulas are defined in  some sublanguages. The following table defines the set of atoms and connectives of the languages of these logics. We call these logics \emph{propositional logics of dependence}.  


{\setlength{\tabcolsep}{0.4em}
\renewcommand{\arraystretch}{1.3}
\begin{center}
\begin{tabular}{|m{5.5cm}|c|c|}
\multicolumn{1}{c}{\textbf{Logic}}&\multicolumn{1}{c}{\textbf{Atoms}}&\multicolumn{1}{c}{\textbf{Connectives}}\\\hline\hline
Propositional dependence logic (\PD)&$p_i,\neg p_i,\bot,\dep(p_{i_1},\dots,p_{i_k},p_j)$&$\wedge,\sor$\\\hline
Propositional dependence logic with intuitionistic disjunction (\PDbor)&$p_i,\neg p_i,\bot$&$\wedge,\sor,\vee$\\\hline
Propositional intuitionistic dependence logic (\PID) &$p_i,\bot,\dep(p_{i_1},\dots,p_{i_k},p_j)$&$\wedge,\vee,\to$\\\hline
Propositional inquisitive logic (\Inql)&$p_i,\bot$&$\wedge,\vee,\to$\\\hline
\end{tabular}
\end{center}
}

\vspace{\baselineskip}

%
%



We consider the above-defined propositional logics of dependence of particular interest for the following reasons. Firstly, \PD and \PID are the propositional logics of their first-order counterparts known in the literature of logics of dependence (see e.g, \cite{Van07dl,Yang2013}).  \Inql is a logic for a new formal semantics, \emph{inquisitive semantics} (\cite{InquiLog}). Although the motivation of \Inql is very different from ours, the logic \Inql adopts essentially the same team semantics. For this reason we regard \Inql as a type of logic of dependence  in this paper. We will discuss the connection between \Inql and \PID in \Cref{sec:inql}. Secondly, we will prove in \Cref{sec:nf} that all of the above logics being syntactically proper fragments of \PT are actually expressively complete as well. We say that two propositional logics of dependence $\LL_1$ and $\LL_2$ have \emph{the same expressive power}, written $\LL_1\equiv\LL_2$, if every formula in the language of $\LL_1$ is semantically equivalent to a formula in the language of $\LL_2$, and vice versa. As a consequence of their expressive completeness, all of the logics mentioned above have the same expressive power. Thirdly, we will show in \Cref{sec:nf} that formulas in the language of these logics have natural disjunctive or conjunctive normal forms that resemble the normal forms of classical propositional logic. On the basis of the normal forms in \Cref{sec:axiom} we axiomatize the logics \PD, \PDbor and \PID. The logic \Inql was already axiomatized in \cite{InquiLog} using an argument that makes essential use of the same normal form. Let us end this section by stating the main results just mentioned as two theorems below.


\begin{theorem}\label{expressive_pw}
All of the logics \PT, \PD, \PDbor, \Inql  and \PID are expressively complete. In particular, $\PT\equiv\PD\equiv\PDbor\equiv\Inql\equiv\PID$.
\end{theorem}

\begin{theorem}
The deduction systems of \PD, \PDbor and \PID defined in \Cref{sec:axiom} are sound and (strongly) complete.
\end{theorem}

\section{Dependence atoms and inquisitive semantics}\label{sec:inql}


In this section, we make some comments on dependence atoms and discuss the connection between inquisitive semantics and logics of dependence.

In general, a (first-order or propositional) dependence atom $\dep(x_1,\dots,x_k,y)$ states the existence of a function $f$ mapping each possible assignment $s$ for $x_{i_1},\dots,x_{i_{k}}$ to an assignment for $y$. In the first-order context, the value $s(x_{i_m})$ ranges over the (possibly infinite) domain of a first-order model. In the propositional context, $s(x_{i_m})$ has only two possible values, namely $0$ or $1$, and the $f$ is a Boolean function that has the finite range $\{0,1\}$. In the presence of the intuitionistic disjunction the behavior of the Boolean function (thus the dependence atom) can be easily expressed. This also explains why we name the logic \PDbor \emph{propositional  dependence logic with intuitionistic disjunction} even though it does not have dependence atoms in its language. In the following lemma we present two equivalent definitions of propositional dependence atoms that are known in the literature (see e.g. \cite{VaMDL08,lovo10,Yang_dissertation,HLSV14}). The proof of their equivalence is left to the reader. Note that the same definitions do not apply to first-order dependence atoms (c.f. \cite{Axiom_fo_d_KV}).



\begin{lemma}\label{dep_v_df}
Let $K=\{i_1,\dots,i_{k}\}$. We have 
\begin{align*}
\dep(p_{i_1},\dots,p_{i_k},p_j)&\equiv \bigsor_{s\in 2^K}\left(p_{i_1}^{s(i_1)}\wedge \dots \wedge p_{i_{k}}^{s(i_{k})}\wedge (p_{j}\bor\neg p_{j})\right)\\
&\equiv \bigbor_{f\in 2^{2^K}}\bigsor_{s\in 2^K}\left(p_{i_1}^{s(i_1)}\wedge \dots \wedge p_{i_{k}}^{s(i_{k})}\wedge p_{j}^{f(s)}\right)
\end{align*}
where $p_i^{1}:=p_i$ and $p_i^0:=\neg p_i$ for any index $i$. In particular,
\(\dep(p_j)\equiv p_j\vee\neg p_j.\)
\end{lemma}

It follows from the above lemma that dependence atoms $\dep(p_{i_1},\dots,p_{i_k},p_j)$ with multiple arguments can be defined using dependence atoms $\dep(p_i)$ with single argument.
We call such  atoms $\dep(p_i)$  \emph{constancy dependence atoms}. The team semantics given in \Cref{TS_PT} for these atoms is reduced to
\begin{itemize}
\item $X\models\,\dep(p_{i})$ iff for all $s,s'\in X$, $s(i)=s'(i)$
\end{itemize}
meaning that $p_i$ has a constant value in the team in question (which can be characterized by the formula $p_i\vee\neg p_i$ as in the preceding lemma). 
An alternative definition of dependence atom with multiple arguments using constancy dependence atoms and intuitionistic implication was given in \cite{AbVan09}. We present this definition without the proof in the following lemma.

 
\begin{lemma}\label{dep_int_df}
$\dep(p_{i_1},\dots,p_{i_k},p_j)\equiv \big(\dep(p_{i_1})\wedge\dots\wedge\dep(p_{i_k})\big)\to\dep(p_{j})$.
\end{lemma}

It follows immediately from the above lemma that propositional intuitionistic dependence logic \PID has the same expressive power as propositional inquisitive logic \Inql. That \Inql independently adopts essentially the same team semantics as propositional logics of dependence was observed by Dick de Jongh and Tadeusz Litak\footnote{In a private conversation with the first author in September 2011.}. It is worthwhile to emphasize that \Inql was introduced with a different motivation along the research line of \emph{inquisitive semantics} (\cite{GR_09,InquiLog}), a new formal semantics characterizing both \emph{assertions} (e.g., \emph{It is raining.}) and \emph{questions} (e.g., \emph{Is it raining?}) in natural language. 
In \Inql a team represents {\em uncertainty} about an assignment. Let us look at an example, the team of \Cref{rer}. The five assignments can be thought to represent uncertainty about one ``right" assignment. In this case there is no uncertainty about $p_1$ because it has value $1$ on every row, but there is some uncertainty about each of $p_2,\dots,p_5$. Subteams represent in this interpretation increased certainty and a singleton team represents absolute certainty. The meaning of the implication $\phi\to\psi$ is here ``Whichever way our knowledge improves so that it supports $\phi$, then also $\psi$ is supported".
The meaning of the dependence atom
$\dep(p_{i_1},\dots,p_{i_k},p_j)$ along the same lines is ``If uncertainty about $p_{i_1},\dots,p_{i_k}$ vanishes, then also $p_j$ becomes certain". The intuition of a team as a state of uncertainty is a useful guiding principle for example in choosing logical operations, and has led to the formulation of $\Inql$ as it is. The intuition in propositional dependence logic is more general covering not only the idea of uncertainty but also the idea of the elements of a team being results of tests or observations without any intention of there being the ``right" assignment.
We refer the reader to \cite{GR_09,InquiLog} for details on inquisitive semantics. 

Strings of the form $\dep(\phi_1,\dots,\phi_k,\psi)$, where the arguments $\phi_1,\dots,\phi_k,\psi$ are arbitrary formulas, are not necessarily well-formed formulas of the logics we study. We decide to apply this restricted syntax for dependence atoms, because we do not have a good intuition for the intended meaning of nested dependence atoms, e.g., strings of the form $\dep(\dep(p_1,p_2),\dep(p_3,p_4))$. A recent work \cite{Ciardelli2015} by Ciardelli interprets in the inquisitive semantics setting the dependence atoms as question entailment, and suggests to define
\[X\models\dep(\phi_1,\dots,\phi_k,\psi)\iff_{\!\!\!\!\textsf{df}}~ X\models\big((\phi_1\vee\neg\phi_1)\wedge\dots\wedge (\phi_{k}\vee\neg\phi_{k})\big)\to (\psi\vee\neg\psi)\]
(which generalizes the equivalences in \Cref{dep_v_df,dep_int_df}).
We do not have a strong argument in favor or against this definition from the perspective of dependence logic, and we are open to other possible alternative definitions for dependence atoms with arbitrary arguments. If we restrict the arguments $\phi_1,\dots,\phi_k,\psi$ to classical formulas, then the above definition is equivalent to
\[X\models\dep(\phi_1,\dots,\phi_k,\psi)\iff_{\!\!\!\!\textsf{df}}~ X\models \bigsor_{s\in 2^K}\left(\phi_{1}^{s(1)}\wedge \dots \wedge \phi_{k}^{s(k)}\wedge (\psi\bor\neg \psi)\right)\]
where $K=\{1,\dots,k\}$.
Dependence atoms with classical arguments are studied in the literature; see e.g., \cite{EHMMVV2013,HLSV14,IemhoffYang15}.

\section{Expressive power and normal forms}\label{sec:nf}

In this section, we study the expressive power and the normal forms of the propositional logics of dependence we defined in Section 1. 

The usual classical propositional logic \CPL is \emph{expressively complete} in the sense that every set $X\subseteq 2^N$ of valuations on $N=\{i_1,\dots,i_n\}$ is definable by some formula $\phi(p_{i_1},\dots,p_{i_n})$ in the language of \CPL, i.e., $X=\{s\in 2^N:s\models\phi\}$. It is well-known that the sets $\{\wedge,\neg\}$, $\{\sor,\neg\}$, $\{\to,\neg\}$ of connectives of \CPL are \emph{functionally complete}. Therefore the fragments of \CPL whose well-formed formulas are defined with these sets of connectives and the atoms $p_i$ are also expressively complete. Analogously, we will show that propositional downward closed team logic \PT and its proper fragments \PD, \PDbor, \PID and \Inql are expressively complete in the sense of \Cref{expres_compl_df} in \Cref{sec:intro}, that is, every nonempty downward closed collection $\mathcal{K}$ of teams  on $N$ is definable by a formula $\phi$ in the language of the logics (i.e., $\mathcal{K}=\{X\subseteq 2^N:X\models\phi\}$). This will prove \Cref{expressive_pw}. The results in  \Cref{expressive_pw} that concern \Inql and \PID can be found essentially in \cite{InquiLog,ivano_msc}; we record them here for the completeness of the paper. The result in  \Cref{expressive_pw} that concerns \PD is due to Taneli Huuskonen who kindly permitted us to include his argument in this paper.

It follows immediately from the expressive completeness result that all of the logics we consider have the same expressive power.
We have  pointed out in the previous section that the expressive equivalence between the two differently motivated logics, \PID and \Inql, is striking, though mathematically this equivalence actually already follows from \Cref{dep_v_df,dep_int_df}. That \PD and \PID are expressively equivalent is also surprising, because their first-order counterparts are very much different from each other in expressive power. First-order dependence logic has the same expressive power as existential second-order logic \cite{Van07dl,KontVan09}, while first-order intuitionistic dependence logic is expressively equivalent to full second-order logic \cite{Yang2013}. 



It will follow from our proof of \Cref{expressive_pw} that every formula in the language of \Inql and \PID is equivalent to a formula in \emph{disjunctive-negative normal form} (which was obtained already in  \cite{InquiLog,ivano_msc}), every formula in the language of \PDbor is equivalent to a formula in  \emph{disjunctive normal form}, 
and every formula in the language of \PD is equivalent to a formula in \emph{conjunctive normal form}.  By \Cref{expressive_pw}, all these disjunctive and conjunctive normal forms  apply to formulas in the language of \PT as well. We therefore will not discuss the normal form for \PT separately. We will argue in the sequel that these normal forms of propositional logics of dependence resemble the same normal forms in classical propositional logic \CPL, because they characterize the truth table of a formula in a similar way. 

We will  treat the logics \PDbor,  \Inql and \PID first, and then the logic \PD. Let us start by defining for each team $X$  a formula $\Theta_X$ in the language of \PDbor and a formula $\Psi_X$ in the language of \Inql and \PID that define the family of subteams of the team $X$. We work with all these three logics at the same time and all results concerning \Inql and \PID can be found essentially in \cite{InquiLog} and \cite{ivano_msc}.


\begin{lemma}[\cite{InquiLog}]\label{X_ThetaX}
Let $X$ be a team on $N=\{i_1,\dots,i_{n}\}$. Define 
\begin{itemize}
\item $\displaystyle\Theta_X:=\bigsor_{s\in X}(p_{i_1}^{s(i_1)}\wedge\dots\wedge p_{i_{n}}^{s(i_{n})})$
\item $\displaystyle\Psi_X:=\neg\neg\bigbor_{s\in X}(p_{i_1}^{s(i_1)}\wedge\dots\wedge p_{i_{n}}^{s(i_{n})})$
\end{itemize}
where we stipulate $\Theta_{\emptyset}=\bot=\Psi_{\emptyset}$.
For any team $Y$ on $N$, we have  \[Y\models \Theta_X\iff Y\subseteq X\iff Y\models \Psi_X.\]
\end{lemma}
\begin{proof}
We first prove that $\Theta_X\equiv\Psi_X$ for all  teams $X$ on $N$. Since $\Theta_X$ is a classical formula and $\Psi_X$ is a negated formula (viewing $\bot$ as $\neg(\bot\to\bot)$), the two formulas are both flat. Thus it suffices to show that $\{s\}\models \Theta_X\iff \{s\}\models \Psi_X$ for all valuations $s$. But this is clear.

It then suffices to check that $Y\models \Theta_X\iff Y\subseteq X$. For ``$\Longrightarrow$'', suppose $Y\models \Theta_X$. If $X=\emptyset$, then $\Theta_X=\bot$, and  $Y=\emptyset=X$. Now assume that $X\neq\emptyset$. For each $s\in X$, there exists a set $Y_s$ such that 
\[Y=\bigcup_{s\in X}Y_s\text{ and }Y_s\models p_{i_1}^{s(i_1)}\wedge\dots\wedge p_{i_{n}}^{s(i_{n})}.\]
Thus either $Y_s=\emptyset$ or $Y_s=\{s\}$ implying $Y\subseteq X$.

Conversely, for ``$\Longleftarrow$'', by the downward closure property it suffices to show that $X\models \Theta_X$. If $X=\emptyset$, then $\Theta_X=\bot$ and $\emptyset\models\bot$. Now, assume that $X\neq\emptyset$. For each $s\in X$, we have  $\{s\}\models p_{i_1}^{s(i_1)}\wedge\dots\wedge p_{i_{n}}^{s(i_{n})}$, which implies that $X\models\Theta_X$.
\end{proof}

We are now ready to give the proof of the claims in \Cref{expressive_pw}  concerning the logics  \PT, \PDbor, \Inql and \PID. The proof for \PD will be given later.
\begin{proof}[Proof of \Cref{expressive_pw} (part 1)]
We will prove that all of the logics \PT, \PDbor, \Inql  and  \PID are expressively complete.
It suffices to show that for every finite set $N=\{i_1,\dots,i_n\}$ of natural numbers, every nonempty downward closed class $\mathcal{K}\in \nabla_N=\{\mathcal{K}\subseteq 2^{2^N}\mid \mathcal{K}\neq \emptyset\text{ and }({X}\in \mathcal{K},~{Y}\subseteq {X}\Longrightarrow {Y}\in\mathcal{K})\}$ of teams on $N$ is definable by a formula $\phi(p_{i_1},\dots,p_{i_n})$  of the logic, i.e., $\mathcal{K}=\llbracket\phi\rrbracket$. 

If $\mathcal{K}=\{\emptyset\}$, then clearly $\mathcal{K}=\llbracket\bot\rrbracket=\llbracket\Theta_{\emptyset}\rrbracket=\llbracket\Psi_{\emptyset}\rrbracket$. Now, assume that $\mathcal{K}\supset \{\emptyset\}$. Since $\mathcal{K}$ is a finite set (of at most $2^{2^n}$ elements), the string $\bigbor_{X\in \mathcal{K}}\Theta_{X}$ is a formula in the language of \PDbor and \PT, and $\bigbor_{X\in \mathcal{K}}\Psi_{X}$ is a formula in the language of \Inql, \PID and \PT. By \Cref{X_ThetaX} for any team $Y$ on $N$ we have 
\[Y\models \bigbor_{X\in \mathcal{K}}\Phi_{X}\iff \exists X\in \mathcal{K}(Y\subseteq X) \iff Y\in \mathcal{K},\]
where $\Phi_X\in \{\Theta_X,\Psi_X\}$. Thus $\llbracket\bigbor_{X\in \mathcal{K}}\Phi_{X}\rrbracket=\mathcal{K}$. 
\end{proof}

 \begin{table}[t]
\centering
\begin{tabular}{c|cc|cccccc}
&$p_{1}$&$p_{2}$&$X_1$&$X_2$&$X_3$&$X_4$&$\cdots$&$X_{16}$\\\hline
$s_1$&$1$&$1$&$0$&$0$&$0$&$0$&$\cdots$&1\\
$s_2$&$0$&$1$&$0$&$1$&$1$&$1$&$\cdots$&1\\
$s_3$&$1$&$0$&$0$&$0$&$1$&$0$&$\cdots$&1\\
$s_4$&$0$&$0$&$0$&$0$&$0$&$1$&$\cdots$&1\\
\hline
$\phi$&&&$1$&$1$&$0$&$1$&$\cdots$&$0$\\
\end{tabular}
\caption{A truth table for a formula $\phi(p_1,p_2)$}
\label{truth_table_example}
\end{table}%

It follows from the proof of the above theorem that every 
formula $\phi(p_{i_1},\dots,p_{i_n})$ 
in the language of \PT is semantically equivalent to a formula in the language of \PDbor in \emph{disjunctive normal form}:
\begin{equation}\label{NF_PDor}
\bigbor_{f\in F}\bigotimes_{s\in X_f}(p_{i_1}^{s(i_1)}\wedge\dots\wedge p_{i_{n}}^{s(i_{n})})
\end{equation}
and also to a formula in the language of \Inql or \PID in \emph{disjunctive-negative normal form}:
\begin{equation}\label{NF_PID}
\bigbor_{f\in F}\neg\neg\bigbor_{s\in X_f}(p_{i_1}^{s(i_1)}\wedge\dots\wedge p_{i_{n}}^{s(i_{n})}),
\end{equation}
where $F$ is a nonempty finite set of indices and each $X_f$ is a  team on $\{i_1,\dots,i_n\}$. It is worthwhile to point out that a similar normal form (of typically infinite size) for first-order dependence logic extended with intuitionistic disjunction was suggested already in \cite{AbVan09}. As with classical propositional logic, the normal form of a formula in the language of \PT or \PDbor or \Inql or \PID can be computed from its truth table. As an illustration, consider the truth table presented in \Cref{truth_table_example} for a formula $\phi(p_1,p_2)$. We first read off  from  the table the set $\mathcal{X}=\{X_f: f\in F\}=\{X: X\models\phi\}$ of teams that satisfy the formula $\phi$. Then for each $f\in F$, we read off the valuations belonging to the team $X_f$ from its corresponding column in the table. Each such column can be viewed as a column in the standard truth table of classical propositional logic, which can be characterized (in the standard way) by a classical formula $\Theta_{X_f}$ in disjunctive normal form. We have, in the end, that $\phi$ is semantically equivalent to a formula $\bigvee_{f\in F}\Theta_{X_f}$ in the language of \PDbor in disjunctive normal form. The set $\mathcal{X}$ in \Cref{truth_table_example} has elements $X_1,X_2,X_4$, etc., and e.g. $\Theta_{X_{4}}=(\neg p_1\wedge p_2)\sor(\neg p_1\wedge\neg p_2)$. In the same manner from the truth table one can also compute the formula $\bigvee_{f\in F}\Psi_{X_f}$  in the language of \Inql or \PID in disjunctive-negative normal form. 




Formulas in disjunctive(-negative) normal form are well-behaved in the sense of the next lemma. In \Cref{sec:axiom} we will apply this lemma to prove the Completeness Theorem for the deduction system of \PDbor.

\begin{lemma}\label{comp_PDbor_main_lm}
Let $\{X_f\mid f\in F\}$, $\{Y_g\mid g\in G\}$ be nonempty finite collections of teams on a set $N$ of natural numbers. Let $\Phi_{X_f}\in\{\Theta_{X_f},\Psi_{X_f}\}$ for each $f\in F$ and $\Phi_{Y_g}\in\{\Theta_{Y_g},\Psi_{Y_g}\}$ for each $g\in G$. The following are equivalent.
\begin{description}
\item[(i)] $\displaystyle\bigbor_{f\in F}\Phi_{X_f}\models\bigbor_{g\in G}\Phi_{Y_g}$.
\item[(ii)] For each $f\in F$, there exists $g_f\in G$ such that $X_f\subseteq Y_{g_f}$.
\end{description}
\end{lemma}
\begin{proof}
%
%
%
(i)$\Longrightarrow$(ii): For each $f_0\in F$ we have  $X_{f_0}\models\Phi_{X_{f_0}}$
by  Lemma \ref{X_ThetaX}. Thus $X_{f_0}\models\bigbor_{f\in F}\Phi_{X_f}$,
which by (i) implies that $ X_{f_0}\models \bigbor_{g\in G}\Phi_{Y_g}$.
This means that
there exists $ g_{f_0}\in G$ such that $X_{f_0}\models \Phi_{Y_{g_{f_0}}}$. Thus by Lemma \ref{X_ThetaX} again we obtain $ X_{f_0}\subseteq Y_{g_{f_0}}$.

(ii)$\Longrightarrow$(i): Suppose $X$ is any team on $N$ satisfying $ X\models\bigbor_{f\in F}\Phi_{X_f}$.
Then $X\models\Phi_{X_f}$ for some $f\in F$, which by \Cref{X_ThetaX} and (ii) means that $X\subseteq X_f\subseteq Y_{g_f}$ for some $ g_f\in G$. Since $Y_{g_f}\models\Phi_{Y_{g_f}}$ holds by Lemma \ref{X_ThetaX}, it follows from the downward closure property that $X\models \Phi_{Y_{g_f}}$, thereby $ X\models \bigbor_{g\in G}\Phi_{Y_g}$, as required.
\end{proof}


Let us now turn to the logic \PD, for which the intuitionistic disjunction is not in the language. In the context of first-order dependence logic the intuitionistic disjunction is uniformly definable using other connectives and quantifiers (see e.g. \cite{Yang_dissertation}). However, this is not the case for propositional dependence logic. The intuitionistic disjunction is not uniformly definable by other connectives in propositional dependence logic \cite{ivano_msc,Yang15}. For this reason we cannot determine the expressive power or obtain the normal form for \PD directly from the results we just established for \PDbor.


We will now give the proof of the claims in \Cref{expressive_pw} concerning \PD and thereby complete the proof of the theorem. Also, a conjunctive normal form for \PD will follow from the proof. This argument is due to Taneli Huuskonen. We present his proof here with his kind permission.

\begin{proof}[Proof of \Cref{expressive_pw} (continued)] 
We will show that \PD is expressively complete.  
It suffices to show that for every finite set $N=\{i_1,\dots,i_n\}$ of natural numbers, every downward closed class $\mathcal{K}\in \nabla_N$ of teams on $N$, there is a formula $\phi(p_{i_1},\dots,p_{i_n})$ in the language of \PD such that $\mathcal{K}=\llbracket \phi\rrbracket$.

Define formulas $\alpha_m$ for each natural number $m$ as follows: 
\[\alpha_0:=\bot,\quad \alpha_1:=\,\dep(p_{i_1})\wedge\dots\wedge\dep(p_{i_n})\quad\text{and}\quad \alpha_m:=\bigsor_{i=1}^m\alpha_1\text{ for }m\geq 2.\]

\begin{claim} 
For any team $X$ on $N$, we have  $X\models\alpha_m$ iff $|X|\leq m$.
\end{claim}

\begin{proofclaim}[Claim 1]Clearly,
\(X\models \alpha_0 \iff X=\emptyset \iff |X|\leq 0\),
and 
\(X\models \alpha_1 \iff  |X|\leq 1.\) If $m>1$, then we have 
\begin{align*}
X\models\alpha_m&\iff X=\bigcup_{i=1}^mX_i\text{ with }X_i\models\alpha_1\text{ for each }i\in\{1,\dots,m\}\\
&\iff X=\bigcup_{i=1}^mX_i\text{ with }|X_i|\leq 1\text{ for each }i\in\{1,\dots,m\}\\
&\iff |X|\leq m.
\end{align*}
\end{proofclaim}


Let $Y$ be a nonempty team on $N$ with $|Y|=k+1$. By Lemma \ref{X_ThetaX} we have 
\begin{equation}\label{Taneli_eq1}
X\models\Theta_{2^N\setminus Y}\iff X\subseteq 2^N\setminus Y\iff X\cap Y=\emptyset.
\end{equation}
Define \(\Xi_Y:=\alpha_k\sor\Theta_{2^N\setminus Y}.\)

\begin{claim} 
For any team $X$ on $N$, we have  $X\models\Xi_Y$ iff $Y\nsubseteq X$.
\end{claim}

\begin{proofclaim}[Claim 2] First, we have 
\begin{align}
X\models\Xi_Y&\iff X=U\cup V\text{ such that }U\models\alpha_k\text{ and }V\models\Theta_{2^N\setminus Y}\nonumber\\
&\iff  X=U\cup V\text{ such that }|U|\leq k\text{ and }V\cap Y=\emptyset. \label{Taneli_eq2}
\end{align}

Next, we prove the claim. Assuming $Y\nsubseteq X$, we have  $|Y\setminus X|\geq 1$ and  $X=(Y\cap X)\cup (X\setminus Y)$, where $(X\setminus Y)\cap Y=\emptyset$. Since
\[|Y\cap X|=|Y|-|Y\setminus X|\leq (k+1)-1=k,\]
by (\ref{Taneli_eq2}) we conclude that $X\models\Xi_Y$.

Conversely, assuming $X\models\Xi_Y$, by (\ref{Taneli_eq2}) we have  $X=U\cup V$, where $|U|\leq k<k+1=|Y|$ and $V\cap Y=\emptyset$. It follows that there exists $s\in Y$ such that $s\notin U\cup V=X$, and therefore $Y\nsubseteq X$.
\end{proofclaim}

Let $\mathcal{K}\in \nabla_N$. Consider the finite class
\(2^{2^N}\setminus \mathcal{K}=\{Y_j\mid j\in J\}\)
of teams on $N$ that are not in $\mathcal{K}$.

\begin{claim} 
For any team $X$ on $N$,
 we have  \(X\in \mathcal{K}\iff Y_j\nsubseteq X\text{ for all }j\in J.\)
\end{claim}

\begin{proofclaim}[Claim 3]
If $X\notin \mathcal{K}$, then $X=Y_{j_0}$ for some $j_0\in J$  by the definition. Thus $Y_{j_0}\subseteq X$. Conversely, if $X\in \mathcal{K}$, then we have  $Y\in \mathcal{K}$ for all $Y\subseteq X$ since  $\mathcal{K}$ is closed under subsets. Hence we must have  $Y_j\nsubseteq X$ for all $j\in J$.
\end{proofclaim}

Finally, since $\emptyset \in \mathcal{K}$, we have  $Y_j\neq \emptyset$ for any $j\in J$. Hence  by Claim 2 and Claim 3 we obtain that for any team $X$ on $N$,
\[X\models \bigwedge_{j\in J}\Xi_{Y_j}\iff Y_j\nsubseteq X\text{ for all }j\in J\iff X\in \mathcal{K},\]
i.e., $\mathcal{K}=\llbracket\bigwedge_{j\in J}\Xi_{Y_j}\rrbracket$, as required.
 \end{proof}
 
We conclude from the proof above that every  formula $\phi(p_{i_1},\dots,p_{i_n})$ in the language of \PD is semantically equivalent to a formula in \emph{conjunctive normal form} 
 \[\bigwedge_{g\in G}\Big(\alpha_{k_g}\sor\Theta_{2^N\setminus Y_g}\Big),\]
where $G$ is a finite set of indices and each $Y_g$ is a nonempty team on $N=\{i_1,\dots,i_n\}$ with $|Y_g|=k_g+1$.  To compute the conjunctive normal form of a formula $\phi$ from its truth table (see \Cref{truth_table_example} for an example), first find the set $\mathcal{Y}=\{Y_g\mid g\in G\}$ of teams that falsify the formula $\phi$. Then we  write the formula $\Theta_{Y_g}^\ast$ according to the column of each $Y_g$  and the formula in normal form is $\bigwedge_{g\in G}\Theta^\ast_{Y_g}$. For example, the set $\mathcal{Y}$ in \Cref{truth_table_example} contains $X_3,X_{16}$, etc., and $\Theta_{X_3}^\ast=(\dep(p_1)\wedge \dep(p_2))\sor(p_1\wedge p_2)\sor (\neg p_1\wedge \neg p_2)$.  


Propositional logic can be viewed as first-order logic over a first-order model $\mathbf{2}$ with a two-element domain $\{0,1\}$. A valuation on an $n$-element index set $N=\{i_1,\dots,i_n\}$ can be viewed as a first-order assignment from a set $\{x_{i_1},\dots,x_{i_n}\}$ of first-order variables into the first-order model $\mathbf{2}$. In this sense \Cref{expressive_pw} shows that for a fixed number $n$, propositional logics of dependence (viewed as a restricted first-order dependence logic) can characterize all nonempty downward closed classes $\mathsf{K}$ of first-order teams of $\mathbf{2}$ with an $n$-element domain. 
These classes $\mathsf{K}$ are called \emph{$(2,n)$-suits} by Cameron and Hodges in \cite{CameronHodges2001}, where the authors estimate  the number of distinct collections $\mathsf{K}$ for fixed $n$. The counting result in  \cite{CameronHodges2001} implies that there is no compositional semantics for the propositional logics of dependence we consider in this paper in which the semantic truth set of a formula $\phi(p_{i_1},\dots,p_{i_n})$ is a subset of $2^N$ justifying that the team semantics given in this paper is indeed an appropriate compositional semantics  for the propositional logics of dependence.

 Another simple but important corollary of \Cref{expressive_pw} is that we can now conclude from the compactness of \Inql, established in \cite{InquiLog}, that all of the propositional logics of dependence we consider are compact.
 
 \begin{theorem}[Compactness]\label{compactness}
The logics \PT, \PD, \PDbor, \Inql  and \PID are compact, that is, if $\Gamma\models\phi$, then there exists a finite subset $\Delta$ of $\Gamma$ such that $\Delta\models\phi$.
\end{theorem}
\begin{proof}
The compactness of \Inql is proved in \cite{InquiLog}. Now, by \Cref{expressive_pw} every formula in the language of any of the other logics is semantically equivalent to a formula in the language of \Inql. Hence the compactness of the other logics follows.
\end{proof}

\section{Axiomatizations}\label{sec:axiom}

We have already noted in Section 1 the significant difference between first-order  dependence logic and propositional logics of dependence. First-order dependence logic  is not axiomatizable as it is equivalent to existential second-order logic (\cite{Van07dl}). Some partial axiomatization of first-order logics of dependence can be found in the literature \cite{Axiom_fo_d_KV,Galliani_general_semantics}. Propositional logics of dependence are axiomatizable in full. This section is devoted to developing complete axiomatizations for the propositional logics of dependence defined in Section 1. The logic \Inql was already axiomatized in \cite{InquiLog}. Adding two obvious axioms characterizing the equivalences in \Cref{dep_v_df,dep_int_df} to the Hilbert style deduction system of \Inql given in  \cite{InquiLog} we obtain  a sound and complete deduction system of \PID  in \Cref{sec:axiom_pid}. In \Cref{sec:axiom_pdv}, we define a natural deduction system of \PDbor and prove the Soundness and Completeness Theorem. Our argument makes essential use of the disjunctive normal form obtained in the previous section. The logic \PD that lacks the intuitionistic disjunction will be axiomatized in \Cref{sec:axiom_pd}.  

A recent paper \cite{SanoVirtema2014} provides a Hilbert style deduction system and a labelled tableau system for \PD and \PDbor. Based on the deduction system in \cite{InquiLog} and the system for \PDbor to be given in \Cref{sec:axiom_pdv} (also recorded in the dissertation \cite{Yang_dissertation} of the first author), in the context of inquisitive semantics a recent paper \cite{Ciardelli2015} gave a sound and complete natural deduction system of the logic \PT without dependence atoms. Adding to this system the obvious rules that characterize the equivalences in \Cref{dep_v_df,dep_int_df}, one obtains immediately a sound and complete natural deduction system of \PT. We will, therefore, not treat \PT in this section.


Before we axiomatize the logics, let us first define and examine the relevant basic notions. A \emph{(syntactic) consequence relation} $\vdash_\LL$ of a logic \LL is a relation on $\mathcal{P}(\textsf{Frm}_\LL)\times \textsf{Frm}_\LL$, where $\textsf{Frm}_\LL$ is the set of all well-formed formulas in the language of  \LL. In what follows we will define Hilbert-style deduction systems or natural deduction systems  for our propositional logics of dependence. For these logics \LL we write $\Gamma\vdash_\LL\phi$ if $\phi\in \Gamma$, or there is a \emph{derivation} (defined as usual) in the given deduction system with the conclusion $\phi$ and the hypotheses from $\Gamma$.  If no confusion arises we write simply $\Gamma\vdash\phi$. If $\vdash_\LL\phi$, then $\phi$ is called a \emph{theorem} of \LL. 
If $\phi\vdash\psi$ and $\psi\vdash\phi$, then we say that $\phi$ and $\psi$ are \emph{provably equivalent}, in symbols $\phi\dashv\vdash\psi$.
The main result of this section is that the  consequence relation $\vdash_\LL$ for every propositional logic of dependence \LL we consider in this paper is sound and (strongly) complete, that is, $\vdash_\LL$ equals $\models$, where $\models$ is the semantic consequence relation defined in Section 1. 

A consequence relation $\vdash_\LL$  is said to be \emph{structural} if $\vdash_\LL$ is closed under uniform substitution, i.e., if $\Gamma\vdash_\LL\phi$, then $\{\sigma(\gamma)\mid \gamma\in \Gamma\}\vdash_\LL\sigma(\phi)$ holds for all substitutions $\sigma$ of \LL. A \emph{substitution} $\sigma$ of  \LL is a function from the set of all propositional variables to the set of all well-formed formulas of \LL that commutes with the connectives in \LL. Formulas of \PD and \PDbor are, by definition, assumed to be in \emph{strict negation normal form}, that is, negation occurs in front of propositional variables only.
This way substitution instances $\neg\sigma(p)$ of the formula $\neg p$ are not necessarily well-formed formulas in the language of \PD or \PDbor. Also, in  \Cref{sec:inql} we mentioned that substitution instances $\dep(\sigma(p_{i_1}),\dots,\sigma(p_{i_k}))$ of the formula $\dep(p_{i_1},\dots,p_{i_k})$ may not be well-formed formulas of the logics either. We already proposed possible ways to extend the logic so as to make the notion of substitution well-defined in our logics. Apart from the definition issue, note that the (sound and complete) consequence relations of these logics cannot be structural. As an illustration, for the logics \PD and \PDbor we have 
\begin{equation}\label{exclmid_exmp}
p_i\sor p_i\models p_i,\text{ whereas }\dep(p_i)\sor\dep(p_i)\not\models\dep(p_i)\text{ and }(p_i\vee\neg p_i)\sor(p_i\vee\neg p_i)\not\models(p_i\vee\neg p_i);
\end{equation}
for the logics \Inql and \PID we have  
\begin{equation}\label{2neg_exmp}
\models \neg\neg p_i\to p_i,\text{ whereas }\not\models \neg\neg (p_i\vee\neg p_i)\to (p_i\vee\neg p_i).
\end{equation}
In other words, the \emph{Substitution Rule} 
\[\AxiomC{$\phi(p_{i_1},\dots,p_{i_n})$}\RightLabel{\Sub} \UnaryInfC{$\phi(\psi_1/p_{i_1},\dots,\psi_n/p_{i_n})$} \DisplayProof\]
is \emph{not} sound in these logics. For this reason in the deduction systems of these logics to be given the axioms or rules should \emph{not} be read as schemata unless otherwise specified. On the other hand, it is shown in \cite{ivano_msc} and \cite{IemhoffYang15} that the consequence relations of these logics are closed under \emph{flat substitutions}, i.e., substitutions $\sigma$ such that $\sigma(p)$ is flat for all propositional variables $p$.

\subsection{\Inql and \PID}\label{sec:axiom_pid}

In this section, we define sound and (strongly) complete Hilbert style deduction systems for propositional intuitionistic dependence logic (\PID) and propositional inquisitive logic (\Inql). All of the results for \Inql were obtained already in \cite{InquiLog}, but we include them in this section for the completeness of the paper and for the benefit of the reader.

As the name suggested, propositional intuitionistic dependence logic \PID satisfies all axioms of propositional intuitionistic logic (\IPC). This is also true for its sublogic propositional inquisitive logic \Inql (which has  the same syntax as \IPC). But as pointed out already, not all the rules of \IPC are sound in \PID and \Inql. Especially the Substitution Rule is not sound in \PID or \Inql. Another property that distinguishes \PID and \Inql from \IPC is that the \emph{(atomic) double negation law} $\neg\neg p_i\leftrightarrow p_i$, of which only the direction $p_i\to\neg\neg p_i$ is valid in \IPC, is valid  in \PID and \Inql. But substitution instances $\neg\neg\phi\leftrightarrow\phi$ of the atomic double negation law are not in general valid in \PID or \Inql; a non-valid instance is given in  (\ref{2neg_exmp}).  It is proven in \cite{ivano_msc} that a formula $\phi$ in the language of \Inql (or \PID) is flat if and only if  $\neg\neg\phi\to \phi$ is valid.


As proved in \cite{InquiLog} and \cite{ivano_msc}, there is a close connection between \Inql and intermediate logics, i.e., logics between intuitionistic propositional logic and classical propositional logic. We briefly discussed in \Cref{sec:intro} that the team semantics induces an intuitionistic Kripke frame. To be more precise, fix a finite set $N$ of natural numbers. The set $\mathcal{P}(2^N)\setminus\{\emptyset\}$ of all nonempty teams on $N$ endowed with the superset relation $\supseteq$ forms a finite Kripke frame $\mathfrak{F}=(\mathcal{P}(2^N)\setminus\{\emptyset\}, \supseteq)$. Define a valuation $V$ on $\mathfrak{F}$ such that a propositional variable $p_i$ is true at a possible world $X$ of $\mathfrak{F}$ if and only if $p_i$ is true on the team $X$. For formulas in the language of \Inql (which is the language of \IPC), the usual single-world-based Kripke semantics corresponds exactly to the team semantics. Moreover, the frame $\mathfrak{F}$ is a so-called \emph{Medvedev frame}, a frame for Medvedev logic  \ML \cite{MedvedevLog}. The \emph{schematic fragment} of \Inql, i.e., the subset of the set of theorems of \Inql that is closed under uniform substitution, is exactly \ML. The logic \Inql itself (the set of theorems) is the \emph{negative variant} of any intermediate logic \LL between Maksimova's logic \ND \cite{MaksimovaLog86} and \ML, such as \ML and Kreisel-Putnam logic \KP \cite{KrsPutnamLog57}. Interested readers are referred to \cite{InquiLog} and  \cite{ivano_msc} for details.

%

Below we present the sound and complete Hilbert style system of \Inql as defined in \cite{InquiLog}. 

\begin{definition}[A Hilbert style deduction system of \Inql, \cite{InquiLog}]\label{Inql_axioms}\
\begin{description}
\item[Axioms:] \
\begin{enumerate}
\item all substitution instances of \IPC axioms
\item   all substitution instances of the Maksimova's axiom $\mathrm{ND}_k$ for all $k\in\mathbb{N}$:
\[(\mathrm{ND}_k)~~~~~~\big(\neg p_j\to\bigvee_{1\leq i\leq k}\neg p_i\big)\to \bigvee_{1\leq i\leq k}(\neg p_j\to\neg p_i).\]
\item $\neg\neg p_i\to p_i$ for all propositional variables $p_i$
\end{enumerate}
\item[Rules:] \
\begin{description}
\item[ Modus Ponens:] \AxiomC{$\phi\to \psi$} \AxiomC{$\psi$}\BinaryInfC{$\psi$} \DisplayProof ~(\MP)
\end{description} 
\end{description}
\end{definition}

Note that the above system does not have the Substitution Rule and the axioms are not axiom schemata. In particular, not every substitution instance of the axiom 3 ($\neg\neg p_i\to p_i$) is sound, whereas the Substitution Rule can be applied to axioms 1 and 2.

One can show in the system of \Inql that every formula $\phi$ in the language of \Inql is provably equivalent to a formula $\bigvee_{f\in F}\Psi_f$ in disjunctive-negative normal form (\ref{NF_PID}). This is a crucial fact used in the proof of the (Weak) Completeness Theorem for the above system given in \cite{InquiLog}. The Strong Completeness Theorem is proved in \cite{ivano_msc} by the standard canonical model argument. The reader is referred to the references for the proof. Below we only state the theorem.

\begin{theorem}[Soundness and Strong Completeness Theorem \cite{InquiLog,ivano_msc}]\label{inql_str_completeness}
For any set $\Gamma\cup\{\phi\}$ of formulas in the language of \Inql, we have $\Gamma\vdash_{\Inql}\phi\iff\Gamma\models\phi$.
\end{theorem}

Adding to the system in \Cref{Inql_axioms} obvious axioms for dependence atoms that correspond to the equivalences in \Cref{dep_v_df,dep_int_df} we obtain a sound and (strongly) complete deduction system of \PID.

\begin{definition}[A Hilbert style deduction system of \PID]\label{PID_axioms}
The system contains all axioms and rules of the system for \Inql given in \Cref{Inql_axioms} together with the following two axioms (not schemata):
\begin{description}
\item[Axioms:] \
\begin{enumerate}
\setcounter{enumi}{3}
\item $\dep(p_i)\leftrightarrow(p_i\vee\neg p_i)$
\item $\dep(p_{i_1},\dots,p_{i_k},p_j)\leftrightarrow \Big(\big(\dep(p_{i_1})\wedge\dots\wedge\dep(p_{i_{k}})\big)\to \dep(p_{j})\Big)$ 
\end{enumerate}
\end{description}
\end{definition}


\begin{theorem}[Soundness and Strong Completeness Theorem]\label{pid_str_completeness}
Let $\Gamma$ be a set of formulas and $\phi$  a formula in the language of \PID. Then $\Gamma\vdash_{\PID}\phi\iff\Gamma\models\phi$.
\end{theorem}
\begin{proof}
Follows from \Cref{inql_str_completeness},  \Cref{dep_v_df,dep_int_df}.
\end{proof}

\subsection{\PDbor}\label{sec:axiom_pdv}

In this section, we define a natural deduction system of propositional dependence logic with intuitionistic disjunction (\PDbor) and prove the Soundness and Completeness Theorem.

Let us first present our natural deduction system of \PDbor. 
As emphasized already, the consequence relation $\vdash_{\PDbor}$ given by this system is \emph{not} structural. In the following definition $\phi$, $\psi$ and $\chi$ are symbols in the metalanguage that stand for arbitrary formulas. 
Axioms and rules presented using concrete formulas such as $p_i$ and dependence atoms should, however, \emph{not} be read as  schemata. 

\begin{definition}[A natural deduction system of \PDbor]\label{Natrual_Deduct_PDbor}\ 

\begin{center}
\setlength{\tabcolsep}{6pt}
\renewcommand{\arraystretch}{1.8}
\setlength{\extrarowheight}{1pt}
\begin{tabular}{|C{11.3cm}|}
\multicolumn{1}{c}{\textbf{AXIOM}}\\\hline\hline
Atomic excluded middle\\
\AxiomC{}\noLine\UnaryInfC{} \RightLabel{\exclmid}\UnaryInfC{$p_i\sor\neg p_i$}\noLine\UnaryInfC{}\DisplayProof\\\hline
\end{tabular}
\end{center}

\begin{center}
\setlength{\tabcolsep}{9pt}
\renewcommand{\arraystretch}{1.8}
\setlength{\extrarowheight}{1pt}
\begin{tabular}{|C{5.4cm}C{5.4cm}|}
\multicolumn{2}{c}{\textbf{RULES}}\\\hline\hline
Conjunction introduction&Conjunction elimination\\
\AxiomC{}\noLine\UnaryInfC{$\phi$} \AxiomC{}\noLine\UnaryInfC{$\psi$}\BinaryInfC{$\phi\wedge\psi$}\noLine\UnaryInfC{} \DisplayProof \ci&\AxiomC{}\noLine\UnaryInfC{$\phi\wedge\psi$} \UnaryInfC{$\phi$}\noLine\UnaryInfC{} \DisplayProof\ce\quad\AxiomC{}\noLine\UnaryInfC{$\phi\wedge\psi$} \UnaryInfC{$\psi$}\noLine\UnaryInfC{} \DisplayProof\ce\\\hline
Intuitionistic disjunction introduction&Intuitionistic disjunction elimination\\
\AxiomC{$\phi$} \UnaryInfC{$\phi\bor\psi$}\noLine\UnaryInfC{}  \DisplayProof\bori\quad
 \AxiomC{$\psi$} \UnaryInfC{$\phi\bor\psi$}\noLine\UnaryInfC{}  \DisplayProof\bori&\AxiomC{$\phi\bor\psi$} \AxiomC{}\noLine\UnaryInfC{$[\phi]$}\noLine\UnaryInfC{$\vdots$}\noLine\UnaryInfC{$\chi$} \AxiomC{$[\psi]$}\noLine\UnaryInfC{$\vdots$}\noLine\UnaryInfC{$\chi$} \RightLabel{\bore}\TrinaryInfC{$\chi$}\noLine\UnaryInfC{} \DisplayProof\\
\hline
Tensor disjunction introduction&Tensor disjunction weak elimination\\
 \AxiomC{}\noLine\UnaryInfC{$\phi$} \UnaryInfC{$\phi\sor\psi$}  \DisplayProof\sori\quad \AxiomC{}\noLine\UnaryInfC{$\phi$} \UnaryInfC{$\phi\sor\psi$}  \DisplayProof\sori&\AxiomC{$\phi\sor\psi$}  \AxiomC{}\noLine\UnaryInfC{$[\phi]$}\noLine\UnaryInfC{$\vdots$}\noLine\UnaryInfC{$\alpha$} \AxiomC{$[\psi]$}\noLine\UnaryInfC{$\vdots$}\noLine\UnaryInfC{$\alpha$} \RightLabel{$(\ast)$~~\sorwe}\TrinaryInfC{$\alpha$}\DisplayProof\\
 &
($\ast$) whenever $\alpha$ is a classical formula\\\hline
Tensor disjunction substitution&Commutative and associative laws for tensor disjunction\\
\multirow{3}{*}{\AxiomC{$\phi\sor\psi$} \AxiomC{$[\psi]$}\noLine\UnaryInfC{$\vdots$}\noLine\UnaryInfC{$\chi$}\noLine\UnaryInfC{} \RightLabel{\sors}\BinaryInfC{$\phi\sor\chi$}\noLine\UnaryInfC{} \DisplayProof}&\AxiomC{}\noLine\UnaryInfC{$\phi\sor\psi$} \UnaryInfC{$\psi\sor\phi$}\noLine\UnaryInfC{} \DisplayProof$\com\sor$\vspace{-24pt}\\
&\\
&\AxiomC{$\phi\sor(\psi\sor \chi)$} \RightLabel{$\ass\sor$}\UnaryInfC{$(\phi\sor\psi)\sor\chi$}\noLine\UnaryInfC{} \DisplayProof\\\hline
Contradiction introduction&Contradiction elimination\\
\AxiomC{}\noLine\UnaryInfC{$p_i\wedge\neg p_i$} \UnaryInfC{$\bot$}\noLine\UnaryInfC{} \DisplayProof\boti&\AxiomC{}\noLine\UnaryInfC{$\phi\sor\bot$}\RightLabel{\bote}\UnaryInfC{$\phi$}\noLine\UnaryInfC{} \DisplayProof\\
\hline
\multicolumn{2}{|c|}{Distributive law}\\
\multicolumn{2}{|c|}{ \AxiomC{}\noLine\UnaryInfC{$\phi\otimes(\psi\bor\chi)$} \UnaryInfC{$(\phi\otimes \psi)\bor(\phi\otimes \chi)$}\noLine\UnaryInfC{}  \DisplayProof$\dstr\sor\bor$} \\\hline
\end{tabular}
\end{center}

\vspace{\baselineskip}
\end{definition}

The above deduction system consists of one axiom and eleven sets of rules. The only axiom, atomic excluded middle \exclmid, is \emph{not} an axiom schema. Especially, as discussed already, the substitution instances of $\neg p_i$ are not necessarily well-formed formulas in the language of  \PDbor. 


The introduction and elimination rules for conjunction $\wedge$ and intuitionistic disjunction $\vee$ are the usual ones. 
Therefore the usual commutative law, associative law and distributive laws, as listed in the next proposition (\Cref{PDbor_derivable_rules}), can be easily derived in the system. The tensor disjunction $\otimes$ behaves classically when applied to classical formulas. So we have the usual introduction and elimination rules for $\otimes$ over classical formulas, but the usual elimination rule is not sound for $\otimes$ in general; a counterexample is given in (\ref{exclmid_exmp}).
Since we only have the weak elimination rule $\sorwe$ for $\sor$,  the substitution, commutative and associative rules for $\sor$ need to be added to the deduction system. We also include a distributive law that involves $\otimes$ in our system, and in \Cref{PDbor_derivable_rules} we derive the other sound distributive laws. Note that not all usual distributive laws that involve $\otimes$ are sound. Among the following three distributive laws
\[(\phi\sor \psi)\wedge(\phi\sor \chi)/\phi\sor(\psi\wedge\chi),\quad (\phi\wedge \psi)\sor(\phi\wedge \chi)/\phi\wedge(\psi\sor\chi)\]
\[\text{ and }(\phi\bor \psi)\sor(\phi\bor \chi)/\phi\bor(\psi\sor\chi)\]
the first two are not sound in general if some of the formulas involved are non-classical and the last one is not sound even for classical formulas. The following clauses are counterexamples for the above three distributive laws with the team $X=\{s_1,s_2\}$, where $s_1(1)=1$ and  $s_2(1)=0$, witnesses their failure:
\[\big((p_1\vee\neg p_1)\sor p_1\big)\wedge \big((p_1\vee\neg p_1)\sor \neg p_1\big)\not\models (p_1\vee\neg p_1)\sor(p_1\wedge\neg p_1)\]
\[\big((p_1\vee\neg p_1)\wedge p_1\big)\sor\big((p_1\vee\neg p_1)\wedge \neg p_1\big)\not\models (p_1\vee\neg p_1)\wedge (p_1\sor\neg p_1)\]
\[(p_1\bor\neg p_1)\sor(p_1\bor\neg p_1)\not\models p_1\bor(\neg p_1\sor\neg p_1)\]

Interesting derivable clauses of our deduction system are listed in the proposition below.

\begin{proposition}\label{PDbor_derivable_rules} The following clauses are derivable:
\begin{description}
\item[(i)]  \textsf{Ex falso}: $\bot\vdash\phi$ 




\item[(ii)]\label{com_ass_c_bor} Commutative  and associative laws for conjunction and intuitionistic disjunction:

\begin{tabular}{ll}
(a) $\phi\wedge\psi\vdash\psi\wedge\phi$ & (b) $\phi\bor\psi\vdash\psi\bor\phi$\vspace{4pt}\\
(c) $(\phi\wedge\psi)\wedge \chi\vdash\phi\wedge(\psi\wedge\chi)$& (d) $(\phi\bor\psi)\bor\chi\vdash\phi\bor(\psi\bor \chi)$
\end{tabular}

\item[(iii)]\label{dstr_bor_c} Distributive laws for intuitionistic disjunction and conjunction:

\begin{tabular}{ll}
(a) $\phi\wedge(\psi\bor\chi)\vdash (\phi\wedge\psi)\bor(\phi\wedge \chi)$& (b) $(\phi\wedge\psi)\bor(\phi\wedge \chi)\vdash\phi\wedge(\psi\bor\chi)$\vspace{4pt}\\
(c) $\phi\bor(\psi\wedge\chi)\vdash(\phi\bor\psi)\wedge(\phi\bor \chi)$& (d) $(\phi\bor\psi)\wedge(\phi\bor \chi)\vdash\phi\bor(\psi\wedge\chi)$
\end{tabular}

%

\item[(iv)]\label{dstr_sor_c} Distributive laws:

\begin{tabular}{l}
(a) $\phi\sor(\psi\wedge\chi)\vdash(\phi\sor \psi)\wedge(\phi\sor \chi)$ ($\dstr\sor\wedge$)\vspace{4pt}\\
(b) $\phi\wedge(\psi\sor\chi)\vdash(\phi\wedge \psi)\sor(\phi\wedge \chi)$ ($\dstr\wedge\sor$)\vspace{4pt}\\
(c) $\phi\bor (\psi\sor \chi)\vdash(\phi\bor\psi)\sor(\phi\bor\chi)$ ($\dstr\bor\sor$)\vspace{4pt}\\
(d) $(\phi\otimes \psi)\bor(\phi\otimes \chi)\vdash\phi\otimes(\psi\bor\chi)$ ($\dstr\sor\bor\sor$)\vspace{4pt}\\
(e) If $\alpha$ is a classical formula, then\vspace{4pt}\\
\quad\begin{tabular}{l}
- $(\alpha\sor \psi)\wedge(\alpha\sor \chi)\vdash\alpha\sor(\psi\wedge\chi)$ ($\dstrs\sor\wedge\sor$)\vspace{4pt}\\
- $(\alpha\wedge \psi)\sor(\alpha\wedge \chi)\vdash\alpha\wedge(\psi\sor\chi)$  ($\dstrs\wedge\sor\wedge$)
\end{tabular}
\end{tabular}
\end{description}
\end{proposition}
\begin{proof}
The items (ii) and (iii) are derived as usual. It remains to derive the other items.

For the item (i), we have the following derivation:
\begin{center}\AxiomC{$\bot$} \RightLabel{\sori}\UnaryInfC{$\phi\sor\bot$}\RightLabel{\bote}\UnaryInfC{$\phi$} \DisplayProof\end{center}


For $\dstr\sor\wedge$, we have the following derivation:
\begin{center}\AxiomC{$\phi\sor(\psi\wedge\chi)$} \AxiomC{[$\psi\wedge\chi$]}\RightLabel{\ce}\UnaryInfC{$\psi$}\RightLabel{\sors}\BinaryInfC{$\phi\sor\psi$}
\AxiomC{$\phi\sor(\psi\wedge\chi)$} \AxiomC{[$\psi\wedge\chi$]}\RightLabel{\ce}\UnaryInfC{$\chi$}\RightLabel{\sors}\BinaryInfC{$\phi\sor\chi$}
\RightLabel{\ci}\BinaryInfC{$(\phi\sor \psi)\wedge(\phi\sor \chi)$} \DisplayProof\end{center}

For $\dstr\wedge\sor$, we have the following derivation:
\begin{center}\AxiomC{[$\psi$]} \AxiomC{$\phi\wedge(\psi\sor\chi)$}\RightLabel{\ce}\UnaryInfC{$\phi$}\RightLabel{\ci}\BinaryInfC{$\phi\wedge\psi$}
\AxiomC{$\phi\wedge(\psi\sor\chi)$}\RightLabel{\ce}\UnaryInfC{$\psi\sor\chi$}
\RightLabel{\sors}\BinaryInfC{$(\phi\wedge \psi)\sor\chi$}
\RightLabel{(by a similar derivation)}\UnaryInfC{$(\phi\wedge \psi)\sor(\phi\wedge\chi)$} \DisplayProof\end{center}

For $\dstr\bor\sor$, we have the following derivation:
\begin{center}
 \def\defaultHypSeparation{\hskip .1in}

\AxiomC{$\phi\bor (\psi\sor \chi)$}

\AxiomC{[$\phi$]}
\RightLabel{\sori}
\UnaryInfC{$\phi\sor\phi$}

\RightLabel{\bori, \sors}
\UnaryInfC{$(\phi\bor\psi)\sor(\phi\bor\chi)$}

\AxiomC{[$\psi\sor\chi$]}

\RightLabel{\bori, \sors}
\UnaryInfC{$(\phi\bor\psi)\sor(\phi\bor\chi)$}

\RightLabel{\bore}
 \TrinaryInfC{$(\phi\bor\psi)\sor(\phi\bor\chi)$}\DisplayProof
\end{center}

For $\dstr\sor\bor\sor$, we have the following derivation:
\begin{center}
\quad\AxiomC{$(\phi\sor\psi)\bor(\phi\sor\chi)$} \AxiomC{$[\phi\sor\psi]$}\AxiomC{$[\psi]$} \RightLabel{\bori}\UnaryInfC{$\psi\bor\chi$}\RightLabel{\sors}\BinaryInfC{$\phi\sor(\psi\bor\chi)$} \AxiomC{$[\phi\sor\chi]$}\AxiomC{$[\chi]$} \RightLabel{\bori}\UnaryInfC{$\psi\bor\chi$}\RightLabel{\sors}\BinaryInfC{$\phi\sor(\psi\bor\chi)$} \RightLabel{\bore}\TrinaryInfC{$\phi\sor(\psi\bor\chi)$} \DisplayProof
\end{center}
 
 For $\dstrs\sor\wedge\sor$, if $\alpha$ is a classical formula, then we have the following derivation:
 \begin{center}

\AxiomC{$(\alpha\sor \psi)\wedge(\alpha\sor \chi)$} 
\RightLabel{$\dstr\wedge\sor$}\UnaryInfC{$\big((\alpha\sor \psi)\wedge\alpha\big)\sor \big((\alpha\sor \psi)\wedge\chi\big)$}

\RightLabel{\ce, \sors}\UnaryInfC{$\alpha\sor \big((\alpha\sor \psi)\wedge\chi\big)$}

\RightLabel{$\dstr\wedge\sor$}\UnaryInfC{$\alpha\sor (\alpha\wedge\chi)\sor (\psi\wedge\chi)$}

\RightLabel{\ce, \sors}\UnaryInfC{$\alpha\sor\alpha\sor (\psi\wedge\chi)$}
\RightLabel{\sorwe, \sors}\UnaryInfC{$\alpha\sor (\psi\wedge\chi)$}
\DisplayProof\end{center}

For $\dstrs\wedge\sor\wedge$, if $\alpha$ is a classical formula, then we have the following derivation:

\begin{center}

\AxiomC{$(\alpha\wedge \psi)\sor(\alpha\wedge \chi)$} 

\RightLabel{\ce, \sors}\UnaryInfC{$\alpha\sor\alpha$}

\RightLabel{\sorwe}\UnaryInfC{$\alpha$}

\AxiomC{$(\alpha\wedge \psi)\sor(\alpha\wedge \chi)$} 

\RightLabel{\ce, \sors}\UnaryInfC{$\psi\sor\chi$}

\RightLabel{\ci}\BinaryInfC{$\alpha\wedge(\psi\sor\chi)$}

\DisplayProof\end{center}
\end{proof}

Next, we prove the Soundness Theorem for our deduction system.

\begin{theorem}[Soundness Theorem]\label{PDbor_soundness}
For any set $\Gamma\cup\{\phi\}$ of formulas in the language of \PDbor, we have 
\(\Gamma\vdash\phi\,\Longrightarrow\,\Gamma\models\phi.\)
\end{theorem}
\begin{proof}
We show that $\Gamma\models\phi$ holds for each derivation  $D=\langle \delta_1,\dots,\delta_k\rangle$ with the conclusion $\phi$ and the hypotheses from $\Gamma$.


If $D=\langle \delta_1\rangle$, then $\phi\in \Gamma$ or $\phi=p_i\sor\neg p_i$. In the former case, obviously $\{\phi\}\models\phi$. The latter case follows from the fact that $\models p_i\sor\neg p_i$.

We only check the inductive step for the rules \sorwe and $\dstr\sor\bor$. The soundness of the other rules can be verified as usual.

\sorwe: Assume that $D_0$, $D_1$ and $D_2$ are derivations for $\Pi_0\vdash\phi\sor\alpha$, $\Pi_1,\phi\vdash\alpha$ and $\Pi_2,\psi\vdash\alpha$, respectively and $\alpha$ is a classical formula. We show that $\Pi_0,\Pi_1,\Pi_2\models\alpha$ follows from the induction hypothesis $\Pi_0\models\phi\sor\psi$, $\Pi_1,\phi\models\alpha$ and $\Pi_2,\psi\models\alpha$.  Suppose $X\models\theta$ for all $\theta\in \Pi_0\cup\Pi_1\cup\Pi_2$. Thus $X\models\phi\sor\psi$, which means that there exist $Y,Z\subseteq X$ such that $X=Y\cup Z$, $Y\models\phi$ and $Z\models\psi$. By the downward closure property we have  $Y\models\theta_1$ for all $\theta_1\in \Pi_1$ and $Z\models\theta_2$ for all $\theta_2\in \Pi_2$. It then follows that $Y\models \alpha$ and $Z\models\alpha$, which yield that $X\models\alpha$, as $\alpha$ is classical (thus flat).

$\dstr\sor\bor$: Assume that $D$ is a derivation for $\Pi\vdash\phi\sor(\psi\bor\chi)$. We will show that $\Pi\models(\phi\sor\psi)\bor(\phi\sor\chi)$ follows from the induction hypothesis $\Pi\models\phi\sor(\psi\bor\chi)$. But this is reduced to showing that $\phi\sor(\psi\bor\chi)\models(\phi\sor\psi)\bor(\phi\sor\chi)$. Now, for any team $X$ such that $X\models\phi\sor(\psi\bor\chi)$, there are teams $Y,Z\subseteq X$ such that $X=Y\cup Z$, $Y\models\phi$ and $Z\models\psi\bor\chi$. It follows that $Y\cup Z\models \phi\sor\psi$ or $Y\cup Z\models \phi\sor\chi$. Hence $X\models(\phi\sor\psi)\bor(\phi\sor\chi)$.
\end{proof}

Clearly, classical propositional logic \CPL is a fragment of \PDbor. We now show that our deduction system of \PDbor is conservative over the usual natural deduction system of \CPL.

\begin{lemma}\label{comp_dep_fr}
For any set $\Gamma\cup\{\phi\}$ of classical formulas,  
\(\Gamma\vdash_\CPL\phi\iff \Gamma\vdash_{\PDbor}\phi\).
\end{lemma}
\begin{proof}
Restricted to classical formulas, our deduction system of \PDbor has the same rules as the usual natural deduction system of \CPC. Thus the direction ``$\Longrightarrow$'' follows. For the converse direction ``$\Longleftarrow$'', assuming $\Gamma\vdash_{\PDbor}\phi$, by the Soundness Theorem we obtain that $\Gamma\models_{\PDbor}\phi$. Thus, for all valuations $s$, $\{s\}\models\gamma$ for all $\gamma\in\Gamma$ implies $\{s\}\models\phi$. But as all the formulas involved are classical formulas, by (\ref{cpl_team_red}) from Section 1, for all valuations $s$, $s\models\gamma$ for all $\gamma\in\Gamma$ implies $s\models\phi$  in the sense of classical propositional logic \CPL. Now, by the Completeness Theorem of \CPL we conclude that $\Gamma\vdash_{\CPL}\phi$.
\end{proof}

We now proceed to prove the main result of this section, the Completeness Theorem for our system. The proof goes through an argument via the normal form that is common in proofs of the Completeness Theorems for logics. 
Already the first proofs of the Completeness Theorem for classical propositional logic by Post \cite{zbMATH02604248} and Bernays \cite{Bernays1926} used the disjunctive normal form. Also, G\"odel's \cite{zbMATH02562387} proof of his Completeness Theorem for classical first order logic was based on an application of a normal form.
In our case, the crucial step for the proof of the Completeness Theorem is to establish that every formula  is provably equivalent in our deduction system to a formula in disjunctive normal form (\ref{NF_PDor}) defined in \Cref{sec:nf}. Let us now state this as a lemma. 

\begin{lemma}\label{DNF_PDbor}
Every formula $\phi(p_{i_1},\dots,p_{i_{n}})$ in the language of \PDbor is provably equivalent to a formula in disjunctive normal form
\begin{equation}\label{NF_PDbor_eq}
\bigbor_{f\in F}\Theta_{X_f},~\text{ where }~\Theta_{X_f}=\bigsor_{s\in X_f}(p_{i_1}^{s(i_1)}\wedge\dots\wedge p_{i_{n}}^{s(i_{n})}),
\end{equation}
 $F$ is a nonempty finite set of indices and each $X_f$ is a team on $\{i_1,\dots,i_{n}\}$. 
\end{lemma}

We will postpone the detailed proof of the above lemma and present the proof of the Completeness Theorem for our system first. 


\begin{theorem}[Completeness Theorem]\label{PDbor_wk_completeness}
For any  formulas $\phi$ and $\psi$  in the language of \PDbor, we have 
\(\psi\models\phi\,\Longrightarrow\,\psi\vdash\phi.\) In particular, $\models\phi\,\Longrightarrow \,\vdash\phi$.
\end{theorem}
\begin{proof}
Suppose $\psi\models\phi$, where $\phi=\phi(p_{i_1},\dots,p_{i_n})$ and $\psi=\psi(p_{i_1},\dots,p_{i_n})$. By Theorem \ref{DNF_PDbor} we have 
\begin{equation}\label{completeness_eq1}
\psi\dashv\vdash \bigbor_{f\in F}\Theta_{X_f}~\text{ and }~\phi\dashv\vdash \bigbor_{g\in G}\Theta_{Y_g}.
\end{equation}
for some nonempty finite sets $\{X_f\mid f\in F\}$ and $\{Y_g\mid g\in G\}$ of teams on $\{i_1,\dots,i_n\}$. The Soundness Theorem implies that
\begin{equation}\label{completeness_eq2}
\bigbor_{f\in F}\Theta_{X_f}\models\bigbor_{g\in G}\Theta_{Y_g}.
\end{equation}

If $\bigbor_{f\in F}\Theta_{X_f}=\Theta_{\emptyset}=\bot$, then we obtain $\psi\vdash\bot\vdash\phi$ by  \exfalso.  If $\bigbor_{g\in G}\Theta_{Y_g}=\Theta_{\emptyset}=\bot$, in view of (\ref{completeness_eq2}) we must have  $\bigbor_{f\in F}\Theta_{X_f}=\Theta_{\emptyset}=\bot$ as well. Hence $\psi\vdash\bot\vdash\phi$ follows from (\ref{completeness_eq1}).


Now assume w.l.o.g. that $\bigbor_{f\in F}\Theta_{X_f},\bigbor_{g\in G}\Theta_{Y_g}\neq\bot$ and $X_f,Y_g\neq \emptyset$ for all $f\in F$ and all $g\in G$. By Lemma \ref{comp_PDbor_main_lm}, for each $f\in F$ we have  $X_f\subseteq Y_{g_f}$ for some $ g_f\in G$. Thus we have the following derivation:
\[\begin{array}{rll}
(1)~ & \Theta_{X_f}&\\
(2)~ & \Theta_{X_f}\sor\Theta_{X_{g_f}\setminus X_f}&\text{(\sori)}\\
(3) ~& \Theta_{Y_{g_f}}&\\
(4) ~& \bigbor_{g\in G}\Theta_{Y_g}& \text{(\bori)}
\end{array}\]
Thus $\Theta_{X_f}\vdash\bigbor_{g\in G}\Theta_{Y_g}$ 
for each $f\in F$, which by the rule \bore implies  $\bigbor_{f\in F}\Theta_{X_f}\vdash\bigbor_{g\in G}\Theta_{Y_g}$. Hence $\psi\vdash\phi$ by (\ref{completeness_eq1}).
\end{proof}

\begin{theorem}[Strong Completeness Theorem]\label{PDbor_completeness}
For any set $\Gamma\cup\{\phi\}$ of formulas in the language of \PDbor, we have 
\(\Gamma\models\phi\,\Longrightarrow\,\Gamma\vdash\phi.\)
\end{theorem}
\begin{proof}
By \Cref{PDbor_wk_completeness} and the Compactness Theorem (\Cref{compactness}).
\end{proof}

For logics that enjoy the Deduction Theorem, the clause $\phi\models\psi\iff\phi\vdash\psi$ in \Cref{PDbor_wk_completeness} is equivalent to the Weak Completeness Theorem (i.e., $\models\phi\iff\vdash\phi$) as $\phi\models\psi$ is equivalent to $\models\phi\to\psi$, and the Strong Completeness Theorem follows readily from the Weak Completeness Theorem and the Compactness Theorem as $\Gamma\models\phi$ is equivalent to $\models\bigwedge\Delta\to \phi$ for some finite subset $\Delta$ of $\Gamma$. These usual arguments, however, do not apply to our logic \PDbor, because \PDbor does not even have an implication in its language (though the Deduction Theorem  with respect to the intuitionistic implication does hold in its extension \PT).  In particular, \Cref{PDbor_wk_completeness}, which is crucial in our proof of the Strong Completeness Theorem via compactness, is actually stronger than the usual Weak Completeness Theorem for the logic. We leave the direct proof of the Strong Completeness Theorem without applying the Compactness Theorem for future work.

 We end this section by supplying the proof of \Cref{DNF_PDbor}. 



\begin{proof}[Proof of \Cref{DNF_PDbor}]
Note that  in the statement of the lemma the fixed set $\{p_{i_1},\dots,p_{i_n}\}$ of variables occurring in the normal form of the formula $\phi$ may be different from the set of variables occurring in the formula $\phi$. In order to take care of this subtle point we first prove the following claim: 

\begin{claim*} If $\{i_1,\dots,i_m\}\subseteq \{j_1,\dots,j_k\}\subseteq \mathbb{N}$, then any formula $\psi(p_{i_1},\dots,p_{i_m})$ in the normal form is provably equivalent to a formula $\theta(p_{j_1},\dots,p_{j_k})$ in the normal form.
\end{claim*}

\begin{proofclaim}
Without loss of generality, we may assume that $K=\{j_1,\dots,j_k\}=\{i_1,\dots,i_m,i_{m+1},\dots,i_k\}$. If $k=m$, then the claim holds trivially. Now assume $k>m$ and
\[\psi(p_{i_1},\dots,p_{i_m})= \bigbor_{f\in F}\bigsor_{s\in X_f}(p_{i_1}^{s(i_1)}\wedge\dots\wedge p_{i_{m}}^{s(i_{m})}),\]
where $F$ is a nonempty finite set of indices and each $X_f$ is a team on $M=\{i_1,\dots,i_m\}$. Let
\[\theta(p_{i_1},\dots,p_{i_k})=\bigbor_{f\in F}\mathop{\bigsor_{s\in 2^K}}_{s\upharpoonright M\in X_f}(p_{i_1}^{s(i_1)}\wedge\dots\dots p_{i_k}^{s(i_k)}).\]

The following derivation proves $\psi\vdash\theta$:
\begin{align*}
(1)~~ & \bigbor_{f\in F}\bigsor_{s\in X_f}(p_{i_1}^{s(i_1)}\wedge\dots\wedge p_{i_{m}}^{s(i_{m})})\\
(2)~~& (p_{i_{m+1}}\sor\neg p_{i_{m+1}})\wedge\dots\wedge (p_{i_k}\sor\neg p_{i_k})~~\text{(\exclmid, \ci)}\\
(3)~~& \bigbor_{f\in F}\Big(\big(\bigsor_{s\in X_f}(p_{i_1}^{s(i_1)}\wedge\dots\wedge p_{i_{m}}^{s(i_{m})})\big)\wedge (p_{i_{m+1}}\sor\neg p_{i_{m+1}})\wedge\dots\wedge (p_{i_k}\sor\neg p_{i_k})\Big)\\
&(\,(1), (2), \ci, \text{ distributive law})\\
(4)~~ & \bigbor_{f\in F}\mathop{\bigsor_{s\in 2^K}}_{s\upharpoonright M\in X_f}(p_{i_1}^{s(i_1)}\wedge\dots\wedge p_{i_{m}}^{s(i_m)}\wedge p_{i_{m+1}}^{s(i_{m+1})}\wedge \dots p_{i_k}^{s(i_k)})~~\text{($\dstr\wedge\sor$, \sori, \sors)}
\end{align*}

Conversely, $\theta\vdash\psi$ is proved by the following derivation:
\begin{align*}
(1)~~ & \bigbor_{f\in F}\mathop{\bigsor_{s\in 2^K}}_{s\upharpoonright M\in X_f}(p_{i_1}^{s(i_1)}\wedge\dots\wedge p_{i_{m}}^{s(i_m)}\wedge p_{i_{m+1}}^{s(i_{m+1})}\wedge \dots p_{i_k}^{s(i_k)})\\
(2)~~ & \bigbor_{f\in F}\mathop{\bigsor_{s\in 2^K}}_{s\upharpoonright M\in X_f}(p_{i_1}^{s(i_1)}\wedge\dots\wedge p_{i_{m}}^{s(i_m)})~~(\ce, \sors)\\
(3)~~ & \bigbor_{f\in F}\bigsor_{s\in X_f}(p_{i_1}^{s(i_1)}\wedge\dots\wedge p_{i_{m}}^{s(i_{m})})~~(\sorwe)
\end{align*}
\end{proofclaim}

Let $N=\{i_1,\dots,i_{n}\}$. We now prove the lemma by induction on $\phi(p_{i_1},\dots,p_{i_n})$.  

\vspace{\baselineskip}

Case $\phi(p_{i_1},\dots,p_{i_n})=p_{i_k}$. Obviously $\displaystyle p_{i_k}\dashv\vdash \bigsor_{s\in \{\mathbf{1}\}}p_{i_k}^{s(i_k)}$, where the function $$\mathbf{1}:\{i_k\}\to \{0,1\}$$  is defined as $\mathbf{1}(i_k)=1$. Then, by the Claim, the formula $p_{i_k}$ is provably equivalent to a formula $\theta(p_{i_1},\dots,p_{i_n})$ in normal form.

\vspace{\baselineskip}

Case $\phi(p_{i_1},\dots,p_{i_n})=\neg p_{i_k}$. This case is proved analogously.

\vspace{\baselineskip}

Case $\phi(p_{i_1},\dots,p_{i_n})=\bot$. Trivially $\bot\dashv\vdash\Theta_{\emptyset}=\bot$.

\vspace{\baselineskip}

 Case $\phi(p_{i_1},\dots,p_{i_n})=\psi(p_{i_1},\dots,p_{i_n})\bor\chi(p_{i_1},\dots,p_{i_n})$. By the induction hypothesis we have 
\begin{equation}\label{PDbor_NF_eq1}
\psi\dashv\vdash\bigbor_{f\in F}\Theta_{X_f}\text{ and }\chi\dashv\vdash\bigbor_{g\in G}\Theta_{X_g},
\end{equation}
where $F,G\neq\emptyset$ and $X_f,X_g\subseteq 2^N$. Applying the rules \bore and \bori we obtain
\[\psi\bor\chi\dashv\vdash\bigbor_{f\in F}\Theta_{X_f}\vee\bigbor_{g\in G}\Theta_{X_g}.\]
If $\psi\dashv\vdash\Theta_\emptyset\dashv\vdash \bot$,  then we obtain further  $\psi\vee\chi\dashv\vdash\bigbor_{g\in G}\Theta_{X_g}$ by  \exfalso and the rules \bore and \bori. Similarly for the case when $\chi\dashv\vdash\Theta_\emptyset$.


\vspace{\baselineskip}

Case $\phi(p_{i_1},\dots,p_{i_n})=\psi(p_{i_1},\dots,p_{i_n})\sor\chi(p_{i_1},\dots,p_{i_n})$. By the induction hypothesis, we have (\ref{PDbor_NF_eq1}).
It suffices to show that $\psi\sor\chi\dashv\vdash \theta$, where
\[\theta= \bigbor_{(f,g)\in F\times G}\,\Theta_{X_f\cup X_g}.\]
 For $\psi\sor\chi\vdash\theta$, we have the following derivation:
\allowdisplaybreaks
\begin{align*}
(1)~~&\psi\sor\chi\\
(2)~~&~\Big(\bigbor_{f\in F}\Theta_{X_f}\Big)\sor\Big(\bigbor_{g\in G}\Theta_{X_g}\Big)\\
(3)~~&\bigbor_{f\in F}\Big(\Theta_{X_f}\sor\Big(\bigbor_{g\in G}\Theta_{X_g}\Big)\Big)
~~\text{($\dstr\sor\bor$)}\\
(4)~~& \bigbor_{f\in F}\bigbor_{g\in G}\Big(\Theta_{X_f}\sor \Theta_{X_g}\Big)
~~\text{($\dstr\sor\bor$)}\\
(5)~~ & \bigbor_{(f,g)\in F\times G}\,\Theta_{X_f\cup X_g}~~\text{(\sorwe, or \bote in case $X_f=\emptyset$ or $X_g=\emptyset$)}
\end{align*}
The other direction $\theta\vdash\psi\sor\chi$ is proved using the rules \sori and $\dstr\sor\bor\sor$.

\vspace{\baselineskip}

Case $\phi(p_{i_1},\dots,p_{i_n})=\psi(p_{i_1},\dots,p_{i_n})\wedge\chi(p_{i_1},\dots,p_{i_n})$. By the induction hypothesis, we have (\ref{PDbor_NF_eq1}), and
\[\Theta_{X_f}=\bigsor_{s_0\in X_f}(p_{i_1}^{s_0(i_1)}\wedge\dots\wedge p_{i_n}^{s_0(i_n)})~\text{ and }~ \Theta_{X_g}=\bigsor_{s_1\in X_g}(p_{i_1}^{s_1(i_1)}\wedge\dots\wedge p_{i_n}^{s_1(i_n)}).\]
If $\psi\dashv\vdash\Theta_\emptyset$ or $\chi\dashv\vdash\Theta_\emptyset$, then we obtain $\psi\wedge\chi\dashv\vdash\bot\dashv\vdash \Theta_\emptyset$ by the rule \ce and \exfalso.
Now, assume that neither of $\psi$ and $\chi$ is equivalent to $\bot$. We show that $\psi\wedge\chi\dashv\vdash \theta$, where
\[\theta= \bigbor_{(f,g)\in F\times G}\,\Theta_{X_f\cap X_g}~ \text{ and }~\Theta_{X_f\cap X_g}=\bigsor_{s\in X_f\cap X_g}(p_{i_1}^{s(i_1)}\wedge\dots\wedge p_{i_n}^{s(i_n)}).\]

For $\psi\wedge\chi\vdash\theta$, we have the following derivation:
\allowdisplaybreaks
\begin{align*}
(1)~~&\psi\wedge\chi\\
(2)~~&\Big(\bigbor_{f\in F}\Theta_{X_f}\Big)\wedge\Big(\bigbor_{g\in G}\Theta_{X_g}\Big)\\
(3)~~ & \bigbor_{f\in F}\bigbor_{g\in G}\Big(\Theta_{X_f}\wedge\Theta_{X_g}\Big)~~~(\text{distributive law})\\
(4)~~ & \bigbor_{f\in F}\bigbor_{g\in G}\bigsor_{s_0\in X_f}\bigsor_{s_1\in X_g}\Big((p_{i_1}^{s_0(i_1)}\wedge\dots\wedge p_{i_n}^{s_0(i_n)})\wedge(p_{i_1}^{s_1(i_1)}\wedge\dots\wedge p_{i_n}^{s_1(i_n)})\Big)\\
&\text{($\dstr\wedge\sor$, \sors)}\\
(5)~~ & \bigbor_{(f,g)\in F\times G}\,\bigsor_{(s_0,s_1)\in X_f\times X_g}\big((p_{i_1}^{s_0(i_1)}\wedge p_{i_1}^{s_1(i_1)})\wedge\dots\wedge( p_{i_n}^{s_0(i_n)}\wedge p_{i_n}^{s_1(i_n)})\big)\\
&\text{(\com\sor, \ass\sor)}\\
(6)~~ & \bigbor_{(f,g)\in F\times G}\Big(\Big(\mathop{\bigsor_{(s_0,s_1)\in X_f\times X_g}}_{s_0=s_1}\big((p_{i_1}^{s_0(i_1)}\wedge p_{i_1}^{s_1(i_1)})\wedge\dots\wedge( p_{i_n}^{s_0(i_n)}\wedge p_{i_n}^{s_1(i_n)})\big)\Big)\\
&\sor\Big(\mathop{\bigsor_{(s_0,s_1)\in X_f\times X_g}}_{s_0\neq s_1}\big((p_{i_1}^{s_0(i_1)}\wedge p_{i_1}^{s_1(i_1)})\wedge\dots\wedge( p_{i_n}^{s_0(i_n)}\wedge p_{i_n}^{s_1(i_n)})\big)\Big)\Big)\\
(7)~~ & \bigbor_{(f,g)\in F\times G}\,\bigsor_{(s,s)\in X_f\times X_g}(p_{i_1}^{s(i_1)}\wedge\dots\wedge p_{i_n}^{s(i_n)})~~(\ce, \bote)\\
(8)~~ & \bigbor_{(f,g)\in F\times G}\,\bigsor_{s\in X_f\cap X_g}(p_{i_1}^{s(i_1)}\wedge\dots\wedge p_{i_n}^{s(i_n)})~~(\sorwe)
\end{align*}

The other direction $\theta\vdash\psi\wedge\chi$ is proved  using the rules \bori, \sori and $\dstrs\wedge\sor\wedge$.
\end{proof}

\subsection{\PD}\label{sec:axiom_pd}\label{arm}









In this section, we develop a sound and complete natural deduction system of propositional dependence logic (\PD). 

In the deduction system to be presented, we will need to specify a particular occurrence of a subformula inside a formula. For this purpose we identify a formula in the language of \PD with a finite string of symbols. A propositional variable $p_i$ is a symbol and the other symbols are $\wedge,\sor,\neg,=,(,)$. We number the symbols in a formula with positive integers starting from the left, as in the following example:
\begin{center}
\begin{tabular}{cccccccccccccccc}
$=$& $($&$p_1$,&$p_2$& $)$&$\sor$&$($&$\neg$&$p_3$&$\wedge$&$\dep$&$($&$p_1$,&$p_2$&$)$&$)$\vspace{4pt}\\
1&2&3&4&5&6&7&8&9&10&11&12&13&14&15&16
\end{tabular}
\end{center}
%
%
If the $m$th symbol of a formula $\phi$ starts a string $\psi$ which is a subformula of $\phi$, we denote the subformula by $[\psi,m]_\phi$, or simply $[\psi,m]$. We will sometimes refer to an occurrence of a formula $\chi$ inside a subformula $\psi$ of $\phi$. In this case, we decide to be sloppy about the notations and will use the same counting for the subformula $\psi$.  We write $\phi(\beta/[\alpha,m])$ for the formula obtained from $\phi$ by replacing the occurrence of the subformula $[\alpha,m]$ by $\beta$. For example, for the  formula $\phi=\dep(p_1,p_2)\sor(\neg p_3\wedge\dep(p_1,p_2))$, we denote the second occurrence of the dependence atom $\dep(p_1,p_2)$ by $[\dep(p_1,p_2),11]$, and the same notation also designates the occurrence of $\dep(p_1,p_2)$ inside the formula $\neg p_3\sor\dep(p_1,p_2)$. The notation $\phi[\beta/[\dep(p_1,p_2),11]]$ designates the formula $\dep(p_1,p_2)\sor(\neg p_3\wedge \beta)$.

Let us now present our natural deduction system of \PD. 
As with the other logics considered in this paper, the consequence relation $\vdash_{\PD}$ of \PD given by the system is \emph{not} structural. In the following definition, only the symbols $\phi$ and $\theta$ stand for arbitrary formulas. The propositional variables $p_i$ in the axioms or rules cannot be substituted by arbitrary formulas.

\begin{definition}[A natural deduction system of \PD]\label{pd_rules} \

\begin{center}
\setlength{\tabcolsep}{14pt}
\renewcommand{\arraystretch}{1.8}
\setlength{\extrarowheight}{2pt}
\begin{tabular}{|C{11cm}|}
 \multicolumn{1}{c}{\textbf{AXIOM}}\\\hline\hline
Atomic excluded middle\\
\AxiomC{} \RightLabel{\exclmid}\UnaryInfC{$p_i\sor\neg p_i$}\noLine\UnaryInfC{}\DisplayProof\\\hline\hline
\end{tabular}
\end{center}

\begin{center}
\setlength{\tabcolsep}{16pt}
\renewcommand{\arraystretch}{1.8}
\setlength{\extrarowheight}{2pt}
\begin{tabular}{|C{11cm}|}
 \multicolumn{1}{c}{\textbf{RULES}}\\\hline\hline
The rules \ci, \ce, \sori, \sorwe, \sors, $\com\sor$, $\ass\sor$, \boti and
 \bote  from \Cref{Natrual_Deduct_PDbor}\\\hline
Dependence atom introduction\\
\AxiomC{$p_i$}\RightLabel{\depiz}\UnaryInfC{$\dep(p_i)$}\DisplayProof\quad\AxiomC{$\neg p_i$}\RightLabel{\depiz}\UnaryInfC{$\dep(p_i)$}\DisplayProof\quad
\def\defaultHypSeparation{\hskip .09in}
 \AxiomC{}\noLine\UnaryInfC{$[\dep(p_{i_1})]$}
 \AxiomC{$\dots$}
 \AxiomC{$[\dep(p_{i_k})]$}
\noLine
\TrinaryInfC{}
\branchDeduce
\DeduceC{$\dep(p_j)$}
\RightLabel{\depik}\UnaryInfC{$\quad\dep( p_{i_1},\dots, p_{i_{k}},p_j)\quad$}\noLine\UnaryInfC{}\DisplayProof
\\\hline
%
Dependence atom elimination\\
\AxiomC{$\dep(p_i)$} \AxiomC{}\noLine\UnaryInfC{$[p_i]$}\noLine\UnaryInfC{$\vdots$}\noLine\UnaryInfC{$\theta$} \AxiomC{$[\neg p_i]$}\noLine\UnaryInfC{$\vdots$}\noLine\UnaryInfC{$\theta$} \RightLabel{\depez}\TrinaryInfC{$\theta$} \DisplayProof\\
 \vspace{-6pt}\\
\def\defaultHypSeparation{\hskip .05in}
\AxiomC{$\dep(p_{i_1},\dots,p_{i_k},p_j)$} \AxiomC{$\dep(p_{i_1})$} \AxiomC{\!\!$\dots$\!\!}\AxiomC{$\dep(p_{i_k})$} \RightLabel{\depek}\QuaternaryInfC{$\dep(p_j)$}\noLine\UnaryInfC{} \DisplayProof\\
\hline
Strong elimination\\
\def\defaultHypSeparation{\hskip .08in}
\AxiomC{$\phi$} \AxiomC{}\noLine\UnaryInfC{[$\phi(p_i/[\dep(p_i),m])$]}\noLine\UnaryInfC{$\vdots$} \noLine\UnaryInfC{$\theta$} \AxiomC{[$\phi(\neg p_i/[\dep(p_i),m])$]}\noLine\UnaryInfC{$\vdots$} \noLine\UnaryInfC{$\theta$} \RightLabel{ \se}\TrinaryInfC{$\theta$}\noLine\UnaryInfC{} \DisplayProof\\\hline
\end{tabular}
\end{center}

\end{definition}

\vspace{\baselineskip}

All rules that do not involve the intuitionistic disjunction $\vee$ in the natural deduction system of \PDbor (\Cref{Natrual_Deduct_PDbor}) are included in the above system of \PD. Thus all clauses in \Cref{PDbor_derivable_rules} that do not involve $\vee$ are derivable in the above system. 

We have discussed in \Cref{sec:inql} that  the behavior of dependence atoms can be characterized completely using intuitionistic connectives: 
\[\dep(p_i)\equiv p_i\bor\neg p_i\text{ and }\dep(p_{i_1},\dots,p_{i_k},p_j)\equiv \big(\dep(p_{i_1})\wedge\dots\wedge\dep(p_{i_{k}})\big)\to \dep(p_{j}).\] In view of the Deduction Theorem, these  equivalences are equivalent to 
\[\dep(p_i)\models p_i\bor\neg p_i,~~~p_i\models\dep(p_i),~~~\neg p_i\models\dep(p_i),\]
\[\big(\dep(p_{i_1})\wedge\dots\wedge\dep(p_{i_{k}})\big)\to \dep(p_{j})\models  \dep(p_{i_1},\dots,p_{i_k},p_j)\] \[\text{and } \dep(p_{i_1},\dots,p_{i_k},p_j), \dep(p_{i_1}),\dots,\dep(p_{i_{k}})\models  \dep(p_{j}).\]
It is not hard to see that the induction and elimination rules for dependence atoms characterize the above clauses. The Strong Elimination rule \se characterizes the interaction between constancy dependence atoms and other formulas, and the rule \depez is a special case of the rule \se. For any formula $\phi$, we have  $\phi\models\phi(p_i\vee\neg p_i/[\dep(p_i),m])$. But as $\vee$ distributes over all connectives, this is reduced to $\phi\models\phi(p_i/[\dep(p_i),m])\vee\phi(\neg p_i/[\dep(p_i),m])$ (see the proof of the Soundness Theorem for details). If the intuitionistic disjunction $\vee$ was a connective in our logic \PD, then we could derive in its sound and complete deduction system using the usual introduction and elimination rules for $\vee$ that $\phi\vdash\phi(p_i/[\dep(p_i),m])\vee\phi(\neg p_i/[\dep(p_i),m])$.
Then, we could derive $\theta$ from $\phi$ as follows: 
 \begin{center}
  \def\ScoreOverhang{0.2pt}
 \AxiomC{$\phi$}
\UnaryInfC{$
\displaystyle
\phi(p_i/[\dep(p_i),m])\vee\phi(\neg p_i/[\dep(p_i),m])$}
 \AxiomC{[$\phi(p_i/[\dep(p_i),m])$]}\noLine\UnaryInfC{$\vdots$} \noLine\UnaryInfC{$\theta$}\AxiomC{[$\phi(\neg p_i/[\dep(p_i),m])$]}\noLine\UnaryInfC{$\vdots$} \noLine\UnaryInfC{$\theta$} 

\TrinaryInfC{$\theta$}
 \DisplayProof \end{center}
Evidently every such derivation can be simulated in our deduction system using the rule \se.



Before we prove the Soundness Theorem, as a demonstration of the system, let us  derive  Armstrong's Axioms \cite{Armstrong_Axioms} in the example below. To simplify notations, we use $x,y,z,\dots$ as metalanguage variables that stand for arbitrary propositional variables.

\begin{example}\label{Armstrong_ax_pd}
The following clauses, known as Armstrong's Axioms, are derivable in the natural deduction system of \PD.
\begin{description}
\item[(i)] $\vdash\dep(x,x)$
\item[(ii)] $\dep(x,y,z)\vdash\dep(y,x,z)$
\item[(iii)] $\dep(x,x,y)\vdash\dep(x,y)$
\item[(iv)] $\dep(y,z)\vdash\dep(x,y,z)$
\item[(v)] $\dep(x,y),\dep(y,z)\vdash\dep(x,z)$
\end{description}
\end{example}
\begin{proof}
The derivations are as follows. For (i): 
\begin{center}
\AxiomC{[$\dep(x)$]}
\RightLabel{\depik}\UnaryInfC{$\dep(x,x)$}\DisplayProof
\end{center}
Note that in the above derivation the formula $\dep(x)$ is both the assumption and the conclusion of a sub-derivation.

For (ii):
\begin{center}
\AxiomC{$\dep(x,y,z)$}\AxiomC{[$\dep(y)$]}\AxiomC{[$\dep(x)$]}\RightLabel{\depek}\TrinaryInfC{$\dep(z)$}\RightLabel{\depik}\UnaryInfC{$\dep(y,x,z)$}\DisplayProof
\end{center}

For (iii):
\begin{center}
\AxiomC{$\dep(x,x,y)$}\AxiomC{[$\dep(x)$]}\AxiomC{[$\dep(x)$]}\RightLabel{\depek}\TrinaryInfC{$\dep(y)$}\RightLabel{\depik}\UnaryInfC{$\dep(x,y)$}\DisplayProof
\end{center}

For (iv):
\begin{center}
\AxiomC{$\dep(y,z)$}\AxiomC{[$\dep(y)$]}\RightLabel{\depek}\BinaryInfC{$\dep(z)$}\RightLabel{\depik}\UnaryInfC{$\dep(x,y,z)$}\DisplayProof
\end{center}

For (v):
\begin{center}
\AxiomC{$\dep(y,z)$}\AxiomC{$\dep(x,y)$}\AxiomC{[$\dep(x)$]}\RightLabel{\depek}\BinaryInfC{$\dep(y)$}\RightLabel{\depek}\BinaryInfC{$\dep(z)$}\RightLabel{\depik}\UnaryInfC{$\dep(x,z)$}\DisplayProof
\end{center}
\end{proof}

Next, we prove the Soundness Theorem for the deduction system of \PD.

\begin{theorem}[Soundness Theorem]\label{PDbor_soundness}
For any  set $\Gamma\cup\{\phi\}$ of formulas in the language of \PD, we have 
\(\Gamma\vdash\phi\,\Longrightarrow\,\Gamma\models\phi.\)
\end{theorem}
\begin{proof}
We show that $\Gamma\models\phi$ holds for each derivation $D$ with conclusion $\phi$ and hypotheses from $\Gamma$.
We only verify the case when the strong elimination rule \se is applied. The other cases are left to the reader.

Put $\phi^\ast_{\mathbf{1}}=\phi(p_i/[\dep(p_i),m])$ and $\phi^\ast_{\mathbf{0}}=\phi(\neg p_i/[\dep(p_i),m])$. Assume that $D_0,D_1$ and $D_2$ are derivations for $\Pi_0\vdash\phi$,\, $\Pi_1,\phi^\ast_{\mathbf{1}}\vdash\theta$ and $\Pi_2,\phi^\ast_{\mathbf{0}}\vdash\theta$, respectively. We show that $\Pi_0,\Pi_1,\Pi_2\models\theta$ follows from the induction hypothesis $\Pi_0\models\phi$,\, $\Pi_1,\phi^\ast_{\mathbf{1}}\models\theta$ and $\Pi_2,\phi^\ast_{\mathbf{0}}\models\theta$. 
This is reduced to showing that $\phi\models\phi^\ast_{\mathbf{1}}\vee\phi^\ast_{\mathbf{0}}$. We prove by induction on subformulas $\psi$ of $\phi$ that $\psi\models\psi^\ast_{\mathbf{1}}\vee\psi^\ast_{\mathbf{0}}$ holds. 

The case when $\psi$ is an occurrence of an atom different from $[\dep(p_i),m]$ trivially holds. Since $\dep(p_i)\models p_i\vee\neg p_i$, the case $\psi=[\dep(p_i),m]$ also holds.

If $\psi=\theta\sor\chi$ and without loss of generality we assume that the occurrence of the dependence atom $[\dep(p_i),m]$ is in the subformula $\theta$, then by the induction hypothesis we have  $\theta\models \theta^\ast_{\mathbf{1}}\vee\theta^\ast_{\mathbf{0}}$. Thus $\theta\sor\chi\models ( \theta^\ast_{\mathbf{1}}\vee\theta^\ast_{\mathbf{0}})\sor\chi\models ( \theta^\ast_{\mathbf{1}}\sor\chi)\vee(\theta^\ast_{\mathbf{0}}\sor\chi)\models\theta^\ast_{\mathbf{1}}\vee\theta^\ast_{\mathbf{0}}$.

The case $\psi=\theta\wedge\chi$ is proved similarly.
\end{proof}

The remainder of this section is devoted to the proof of the Completeness Theorem for our system. Our argument makes heavy use of the notion of the \emph{realization} of a formula. We will give the definition of realizations in two steps. Firstly we define realizations for dependence atoms. 

Let $\alpha=\dep(p_{i_1},\dots,p_{i_k},p_j)$ be a dependence atom and $K=\{i_1,\dots,i_k\}$. We call a function $f:2^K\to 2$ a \emph{realizing function} for $\alpha$, where we stipulate that $2^\emptyset=\{\emptyset\}$.  A \emph{realization} $\alpha^\ast_f$ of the dependence atom $\alpha$ over $f$ is defined as
\[\alpha_{f}^\ast:=\displaystyle\bigsor_{s\in 2^K}\left(p_{i_1}^{s(i_1)}\wedge \dots \wedge p_{i_k}^{s(i_k)}\wedge p_j^{f(s)}\right).\]
By \Cref{dep_v_df}, $\alpha$ is semantically equivalent to the intuitionistic disjunction of all its realizations, i.e., $\alpha\equiv\bigvee_{f\in 2^{2^K}}\alpha_f^\ast$. The function $f$ in a realization $\alpha_f^\ast$ characterizes one type of functional dependency between $p_j$ and $p_{i_1},\dots,p_{i_k}$, and thus the formula $\alpha^\ast_f$ ``realizes'' such  dependency. Here the formula $\bigvee_{f\in 2^{2^K}}\alpha_f^\ast$  is not in the language of \PD, as the intuitionistic disjunction $\vee$ is not an eligible connective in the logic \PD. Nevertheless, we will show that in the deduction system of \PD one derives \emph{essentially}, in a sense that becomes apparent below,   the equivalence between $\alpha$ and $\bigvee_{f\in 2^{2^K}}\alpha_f^\ast$. We first prove the easy direction that $\alpha$ essentially follows from $\bigvee_{f\in 2^{2^K}}\alpha_f^\ast$. This, in our deduction system of \PD, is simulated by the derivation that $\alpha$ follows from each intuitionistic disjunct $\alpha_f^\ast$. 



\begin{lemma}\label{depezk_derive}
 If $f$ is a realizing function for a dependence atom $\alpha$, then $\alpha_f^\ast\vdash \alpha$.
\end{lemma}
\begin{proof}
Suppose that $f$ is a realizing function for $\alpha=\dep(p_{i_1},\dots,p_{i_k},p_j)$. We will show 
\begin{equation}\label{depezk_derive_eq1}
\alpha_f^\ast,\dep(p_{i_1}),\dots, \dep(p_{i_k})\vdash \dep(p_j),
\end{equation}
from which $\alpha_f^\ast\vdash \alpha$ will follow by applying the rule \depik.

\begin{claim*} 
For all $t:\{1,\dots,k\}\to 2$, we have  $\alpha_f^\ast,p_{i_1}^{t(1)},\dots,p_{i_k}^{t(k)}\vdash \dep(p_j)$.
\end{claim*}
\begin{proofclaim}
If $i_a=i_b$ and $t(a)=1\neq 0=t(b)$ for some $a,b\in\{1,\dots, k\}$, then 
\[\alpha_f^\ast,p_{i_1}^{t(1)},\dots,p_{i_k}^{t(k)}\vdash p_{i_a}\wedge \neg p_{i_b}\vdash\bot\vdash\dep(p_j)\] 
by the rule \boti and \exfalso. Now, assume $t(a)= t(b)$ for all $a,b\in\{1,\dots, k\}$ with $i_a=i_b$. The derivation 
is as follows:
\[\begin{array}{rll}
(1)& p_{i_1}^{t(1)}\wedge\dots\wedge p_{i_k}^{t(k)}&\text{(assumption, \ci)}\\
(2)& \displaystyle\bigsor_{s\in 2^K}\left(p_{i_1}^{s(i_1)}\wedge \dots \wedge p_{i_k}^{s(i_k)}\wedge p_j^{f(s)}\right)&\text{(assumption)}\\
(3)& \displaystyle\bigsor_{s\in 2^K}\left((p_{i_1}^{t(1)}\wedge\dots\wedge p_{i_k}^{t(k)})\wedge p_{i_1}^{s(i_1)}\wedge \dots \wedge p_{i_k}^{s(i_k)}\wedge p_j^{f(s)}\right)&\text{((1), (2), \ci, $\dstr\wedge\sor$)}\\
(4)& \displaystyle (p_{i_1}^{t(1)}\wedge\dots\wedge p_{i_k}^{t(k)}\wedge p_j^{f(s_t)})\sor \bigsor_{s\in 2^K\setminus\{s_t\}}\bot,&\\
&\text{ where $s_t(i_n)=t(n)$ for each $n\in\{1,\dots, k\}$}&\text{(\boti, \ce, \sors)}\\
(5)& p_j^{f(s_t)}&\text{(\bote, \ce)}\\
(6)&\dep(p_j)&\text{(\depiz)}
\end{array}\]
\end{proofclaim}

Next, we prove (\ref{depezk_derive_eq1}). By the Claim, for all $s_1:\{2,\dots,k\}\to 2$ we have 
\[\alpha_f^\ast,p_{i_1},p_{i_2}^{s_1(2)},\dots,p_{i_k}^{s_1(k)}\vdash\dep(p_j)~~\text{ and }~~\alpha_f^\ast,\neg p_{i_1},p_{i_2}^{s_1(2)},\dots,p_{i_k}^{s_1(k)}\vdash\dep(p_j).\] 
Then, by the rule \depez, we derive that for all $s_1:\{2,\dots,k\}\to 2$,
\begin{equation*}
\alpha_f^\ast,\dep(p_{i_1}),p_{i_2}^{s_1(2)},\dots,p_{i_k}^{s_1(k)}\vdash\theta
\end{equation*}
Repeating this argument $k$ times, we obtain (\ref{depezk_derive_eq1}). This completes the proof.
\end{proof}

The converse direction that $\alpha$ essentially implies  $\bigvee_{f\in 2^{2^K}}\alpha_f^\ast$ will be a special case of a  general result (\Cref{sek_derive}) in what follows. Now, we are ready to define realizations for arbitrary formulas $\phi$.
Let 
\[o=\langle [\alpha_1,m_1],\,\dots,\,[\alpha_c,m_c]\rangle\] 
be a sequence of some occurrences of dependence atoms in  $\phi$.
A \emph{realizing sequence of $\phi$ over $o$} is a sequence $\Omega=\langle f_1,\dots,f_{c}\rangle$ such that each $f_i$ is a realizing function for $\alpha_i$. We call the formula $\phi_\Omega^\ast$ defined 
as follows a \emph{realization of $\phi$ over $o$}: 
\[\phi_{\langle f_1,\dots,f_{c}\rangle}^\ast:=\phi((\alpha_1)_{f_1}^\ast/[\alpha_1,m_1],\dots,(\alpha_c)_{f_c}^\ast/[\alpha_c,m_c]).\]
Let $O$ be the sequence of all of the occurrences of dependence atoms in $\phi$. A realizing sequence of $\phi$ over $O$ is called a \emph{maximal realizing sequence}, and a realization $\phi^\ast_\Omega$ over $O$ is called a \emph{complete realization} of $\phi$. Note that a complete realization $\phi^\ast_\Omega$ of a formula $\phi$ in the language of \PD is a classical formula. For example, consider the formula $\phi=\dep(p_1,p_2)\sor(\neg p_3\wedge\dep(p_1,p_2))$ that we discussed earlier. Let $\epsilon,\delta:2^{\{1\}}\to 2$ be two realizing functions for $\dep(p_1,p_2)$. Over $o=\langle [\dep(p_1,p_2),11]\rangle$ the sequence $\langle \epsilon\rangle$ is a realizing sequence of $\phi$ and the formula $\dep(p_1,p_2)\sor(\neg p_3\wedge(\dep(p_1,p_2))^\ast_\epsilon)$  is a realization of $\phi$. Both $(\dep(p_1,p_2))^\ast_\delta\sor(\neg p_3\wedge(\dep(p_1,p_2))^\ast_\epsilon)$ and $(\dep(p_1,p_2))^\ast_\epsilon\sor(\neg p_3\wedge(\dep(p_1,p_2))^\ast_\delta)$ are complete realizations of $\phi$. Note that if a formula does not contain dependence atoms, then its maximal realizing sequence is the (unique) empty sequence $\langle\rangle$.


Having defined the crucial notion for the proof of the Completeness Theorem, we now check in the next lemma that every formula is semantically equivalent to the intuitionistic disjunction of  its realizations over some $o$. In the proof of our main theorem (\Cref{PD_completeness}), we are only interested in complete realizations, but the  lemmas leading to this theorem will be stated in their most general forms.

\begin{lemma}\label{WNF_PD_sound}
Let $\phi$ be a formula in the language of \PD and $\Lambda$ the set of all its realizing sequences over a sequence $o$. 
Then 
\(\displaystyle\phi\equiv \bigbor_{\Omega\in \Lambda}\phi^\ast_\Omega.\)
\end{lemma}
\begin{proof}
We prove the lemma by induction on $\phi$. 
The case that $\phi$ is a dependence atom follows from \Cref{dep_v_df}. The case that $\phi$ is some other atom is trivial.


The induction case $\phi=\psi\sor\chi$ follows from the induction hypothesis and the fact that 
\([\,A\models A'\text{ and }B\models B'\,] \Longrightarrow A\sor B\models A'\sor B'\)
and that
$A\sor(B\bor C)\models (A\sor B)\bor(A\sor C)$. Similarly for the case $\phi=\psi\wedge\chi$.
\end{proof}

A crucial fact to be used in the proof of the Completeness Theorem is that  a formula $\phi$ is essentially provably equivalent to  $\bigbor_{\Omega\in \Lambda}\phi^\ast_\Omega$. We first prove the easy direction that $\phi$ follows essentially from $ \bigbor_{\Omega\in \Lambda}\phi^\ast_\Omega$. As with realizations for dependence atoms, we will show, in the absence of the intuitionistic disjunction, that each  realization of $\phi$ implies $\phi$.

\begin{lemma}\label{si_derive} 
If $\Omega$ is a realizing sequence of a formula $\phi$ over a sequence of some occurrences of dependence atoms in $\phi$, then $\phi^\ast_\Omega\vdash\phi$.
\end{lemma}
\begin{proof}
We prove the lemma by induction on $\phi$. 
The case that $\phi$ is a dependence atom follows from \Cref{depezk_derive}. The case that $\phi$ is some other atom is trivial.


The induction case $\phi=\psi\sor\chi$ is proved by applying the rule \sors and the case $\phi=\psi\wedge\chi$  is proved by applying  the rule \cs.
\end{proof}

We now turn to the converse direction  that $\phi$ essentially  implies $\bigbor_{\Omega\in \Lambda}\phi^\ast_\Omega$. Let us ponder on what we should settle here in the logic \PD in the absence of the intuitionistic disjunction. Assume that the intuitionistic disjunction was a connective in our logic and we have proved that $\phi\vdash \bigbor_{\Omega\in \Lambda}\phi^\ast_\Omega$. If we have a derivation of $\phi$, then to derive a formula $\theta$ from the same assumptions, it is sufficient to prove that $\theta$ follows from each realization $\phi^\ast_\Omega$ of $\phi$. Our idea is to discard the intermediate step that involves the intuitionistic disjunction and conclude $\theta$ directly from $\phi$ as soon as the fact  that $\theta$ follows from each $\phi^\ast_\Omega$ is settled. To be precise, we will prove the lemma below.

\begin{lemma}\label{sek_derive}\
Let $\Lambda$ be the set of all realizing sequences of $\phi$ over a sequence $o$ of some occurrences of dependence atoms in a formula $\phi$. Let $\{\Pi_\Omega\mid \Omega\in\Lambda\}$ be a collection of sets of formulas. For any formula $\theta$, if 
\begin{equation}\label{sek_derive_eq1}
\forall \Omega\in \Lambda\,(\Pi_\Omega,\phi^\ast_{\Omega}\vdash\theta),
\end{equation}
then $\bigcup_{\Omega\in \Lambda}\Pi_\Omega,\phi\vdash\theta$.
\end{lemma}

We will now work towards a proof of the above lemma. Clearly, if we consider only a single occurrence of a constancy dependence atom in a formula, i.e., if $o=\langle [\dep(p_i),m]\rangle$, then the implication in the lemma follows immediately by applying the rule \se. The idea is then to reduce the general case to this special case in a certain way.

First we show that dependence atoms with multiple arguments can be reduced to  formulas with constancy dependence atoms only.
If $\alpha=\dep(p_{i_1},\dots,p_{i_k},p_j)$, we write
\[\displaystyle\alpha^\ast:=\bigsor_{s\in 2^K}\left(p_{i_1}^{s(i_1)}\wedge \dots \wedge p_{i_k}^{s(i_k)}\wedge \dep(p_j)\right),\]
where $K=\{i_1,\dots,i_k\}$. By \Cref{dep_v_df}, $\alpha$ and $\alpha^\ast$ are semantically equivalent.
Let $o=\langle [\alpha_1,m_1],\,\dots,\,[\alpha_c,m_c]\rangle$ be a sequence of some occurrences of dependence atoms in a formula $\phi$ and let
\[\phi^\ast_o:=\phi(\alpha_1^\ast/[\alpha_1,m_1],\dots,\alpha_c^\ast/[\alpha_c,m_c]).\]


\begin{lemma}\label{phi_derive_phiast}
Let $\phi$ and $o$ be as above. Then $\phi\vdash\phi^\ast_o$.
\end{lemma}
\begin{proof}
We prove the lemma by induction on $\phi$. If $\phi=\dep(p_{i_1},\dots,p_{i_k},p_j)$ and $K=\{i_1,\dots,i_k\}$, then we first derive 
\begin{equation}\label{pd_nf_derive_eq3}
\dep(p_{i_1},\dots,p_{i_k},p_j)\vdash\bigsor_{s\in 2^K}\left(\dep(p_{i_1},\dots,p_{i_k},p_j)\wedge p_{i_1}^{s(i_1)}\wedge \dots \wedge p_{i_k}^{s(i_k)}\right)
\end{equation}
as follows:
\[\begin{array}{rll}
(1)& \dep(p_{i_1},\dots,p_{i_k},p_j)&\text{(assumption)}\\
(2)&(p_{i_1}\sor \neg p_{i_1})\wedge \dots\wedge (p_{i_k}\sor\neg p_{i_k})&\text{(\exclmid, \ci)}\\
(3)&\displaystyle\bigsor_{s\in 2^K}\left(p_{i_1}^{s(i_1)}\wedge \dots \wedge p_{i_k}^{s(i_k)}\right)&\text{($\dstr\wedge\sor$)}\\
(4)&\displaystyle\bigsor_{s\in 2^K}\left( \dep(p_{i_1},\dots,p_{i_k},p_j)\wedge p_{i_1}^{s(i_1)}\wedge \dots \wedge p_{i_k}^{s(i_k)}\right)&\text{((1), (3), $\dstr\wedge\sor$)}\\
\end{array}\]
Next, for each $s:K\to 2$, by the rules \depiz and \depek, we have 
\begin{align*}
\dep(p_{i_1},\dots,p_{i_k},p_j)\wedge p_{i_1}^{s(i_1)}\wedge \dots \wedge p_{i_k}^{s(i_k)}&\vdash\dep(p_{i_1},\dots,p_{i_k},p_j)\wedge \dep(p_{i_1})\wedge\dots\wedge \dep(p_{i_k})\\
&\vdash \dep(p_j)
\end{align*}
Hence $\dep(p_{i_1},\dots,p_{i_k},p_j)\vdash(\dep(p_{i_1},\dots,p_{i_k},p_j))^\ast_o$ follows from (\ref{pd_nf_derive_eq3}) and \sors.


If $\phi\in\{p_i,\neg p_i,\bot\}$, then the statement is trivial. The induction steps can be derived using the rules \ci and \sors.
\end{proof}

In the next lemma we derive the special case of \Cref{sek_derive} in which the sequence $o$ consists of occurrences of constancy dependence atoms only.
For a constancy dependence atom $\dep(p_j)$ there are only two realizing functions, namely the constant functions $\mathbf{1}:\{\emptyset\}\to 2$ and $\mathbf{0}:\{\emptyset\}\to 2$ defined as $\mathbf{1}(\emptyset)=1$ and $\mathbf{0}(\emptyset)=0$. And $\dep(p_j)_{\mathbf{1}}^\ast=p_j$ and $\dep(p_j)_{\mathbf{0}}^\ast=\neg p_j$. A realizing sequence $\Omega=\langle f_1,\dots,f_k\rangle$ over a sequence $o$ of occurrences of constancy dependence atoms is a sequence of constant functions $\mathbf{1}$ and $\mathbf{0}$, and  can, therefore, be characterized by a function $g_\Omega:\{1,\dots,k\}\to 2$ defined as $g_{\Omega}(i)=f_i(\emptyset)$. We identify the sequence $\Omega$ with the function $g_\Omega$ and sometimes write $\phi^\ast_{g_\Omega}$ instead of $\phi^\ast_{\Omega}$. 

\begin{lemma}\label{sezk_derive}
For each function $g: \{1,\dots,k\}\to 2$, Let $\Pi_g$ be a set of formulas. For any formulas $\phi$ and $\theta$,
if  $\Pi_g,\phi^\ast_g\vdash\theta$ holds for all $g: \{1,\dots,k\}\to 2$, i.e., 
\[\Pi_g,\phi(p_{j_1}^{g(1)}/[\dep(p_{j_1}),m_1],\dots, p_{j_k}^{g(k)}/[\dep(p_{j_k}),m_k])\vdash\theta,\] 
then $\bigcup_{g}\Pi_g,\phi\vdash\theta$.
\end{lemma}
\begin{proof}
By the assumption,  for all $h_1:\{2,\dots,k\}\to 2$ we have 
\[\Pi_{g_{1,h_1}},\phi(p_{j_1},p_{j_2}^{h_1(2)}/[\dep(p_{j_2}),m_2],\dots, p_{j_k}^{h_1(k)}/[\dep(p_{j_k}),m_k])\vdash\theta\] 
and 
\[\Pi_{g_{0,h_1}},\phi(\neg p_{j_1},p_{j_2}^{h_1(2)}/[\dep(p_{j_2}),m_2],\dots, p_{j_k}^{h_1(k)}/[\dep(p_{j_k}),m_k])\vdash\theta,\]
where $g_{1,h_1}:\{1,\dots,k\}\to 2$ is the function $h_1\cup\{(1,1)\}$ and $g_{0,h_1}:\{1,\dots,k\}\to 2$ is the function $h_1\cup\{(1,0)\}$.
Then by the rule \se we derive  for all $h_1:\{2,\dots,k\}\to 2$,
\begin{equation*}
\Pi_{g_{1,h_1}},\Pi_{g_{0,h_1}},\phi(p_{j_2}^{h_1(2)}/[\dep(p_{j_2}),m_2],\dots, p_{j_k}^{h_1(k)}/[\dep(p_{j_k}),m_k])\vdash\theta.
\end{equation*}
Repeating this argument $k$ times, we obtain $\bigcup_g\Pi_g,\phi\vdash\theta$ as required.
\end{proof}

Having analyzed the simple realizations that involve constancy dependence atoms only, let us now give the proof of \Cref{sek_derive}. Our main argument is to apply \Cref{phi_derive_phiast} to transform dependence atoms with multiple arguments to formulas with constancy dependence atoms only so as to apply \Cref{sezk_derive} to finish the proof.

\begin{proof}[Proof of \Cref{sek_derive}]
Let $o=\langle[\alpha_1,m_1],\dots,[\alpha_c,m_c] \rangle$,  
\(\alpha_i=\dep(p_{i_1},\dots,p_{i_{k_i}},p_{j_i})\) and $K_i=\{i_1,\dots,i_{k_i}\}$ for any $i\in\{1,\dots,c\}$.
In view of \Cref{phi_derive_phiast}, to prove the lemma it suffices to show that the assumption (\ref{sek_derive_eq1}) implies $\bigcup_{\Omega\in \Lambda}\Pi_\Omega,\phi^\ast_o\vdash\theta$.

For each $i\in\{1,\dots ,c\}$, let 
\(o_i=\langle [\dep(p_{j_i}),l_{i,s}] \rangle_{s\in 2^{K_i}}\)
be the sequence of all occurrences of the (constancy) dependence atoms in the subformula 
\[\alpha^\ast_i=\bigsor_{s\in 2^{K_i}}\left(p_{i_1}^{s(i_1)}\wedge \dots \wedge p_{i_{k_i}}^{s(i_{k_i})}\wedge \dep(p_{j_i})\right)\]
of $\phi^\ast_o$. Each function $g:\bigcup_{i=1}^c\{(i,s)\mid s\in 2^{K_i}\}\to 2$ induces a realizing sequence of $\phi^\ast_o$ over the sequence $\langle o_1,\dots,o_c\rangle$ of occurrences of constancy dependence atoms. For each $i\in \{1,\dots,c\}$, define a function $f_{g,i}:2^{K_i}\to 2$ as
\[f_{g,i}(s)=g((i,s)).\]
Observe that $f_{g,i}$ is a realizing function for $\alpha_i$ and $\Omega_g=\langle f_{g,1},\dots,f_{g,c}\rangle$ is a realizing sequence of $\phi$ over $o$.

\begin{claim*}
$\Pi_{\Omega_g},(\phi^\ast_o)^\ast_{g}\vdash\theta$.
\end{claim*}
\begin{proofclaim}
For each $i\in\{1,\dots,c\}$, we have  $g\upharpoonright \{(i,s)\mid s\in 2^{K_i}\}$ induces a realizing sequence of $\alpha_i^\ast$ over $o_i$ and 
\[(\alpha_i)_{f_{g,i}}^\ast=\displaystyle\bigsor_{s\in 2^{K_i}}\left(p_{i_1}^{s(i_1)}\wedge \dots \wedge p_{i_{k_i}}^{s(i_{k_i})}\wedge p_{j_i}^{f_{g,i}(s)}\right)=(\alpha_i^\ast)^\ast_{g}.\]
Now, by (\ref{sek_derive_eq1}) we  have  
\[\Pi_{\Omega_g},\phi((\alpha_1)_{f_{g,1}}^\ast/[\alpha_1,m_1],\dots,(\alpha_c)_{f_{g,c}}^\ast/[\alpha_c,m_c])\vdash\theta,\] 
which is
\begin{equation*}
\Pi_{\Omega_g},\phi((\alpha_1^\ast)^\ast_{g}/[\alpha_1,m_1],\dots,(\alpha_c^\ast)_{g}^\ast/[\alpha_c,m_c])\vdash\theta,\text{ i.e., }\Pi_{\Omega_g},(\phi^\ast_o)^\ast_{g}\vdash\theta.
\end{equation*}
\end{proofclaim}

Put $\bigcup_{i=1}^c\{(i,s)\mid s\in 2^{K_i}\}=A$. By the Claim and \Cref{sezk_derive}, we obtain that $\bigcup_{g\in 2^A}\Pi_{\Omega_g},\phi^\ast_o\vdash\theta$. It is not hard to see that there is a one-to-one correspondence between functions $g:A\to 2$ and realizing sequences $\Omega$ of $\phi$ over $o$. This means that $\bigcup_{g\in 2^A}\Pi_{\Omega_g}=\bigcup_{\Omega\in \Lambda}\Pi_\Omega$, which completes the proof.
\end{proof}

Now, we are in a position to prove the (Strong) Completeness Theorem for the deduction system of \PD.

\begin{theorem}[Strong Completeness Theorem]\label{PD_completeness}
For any  set $\Gamma\cup\{\phi\}$ of formulas in the language of \PD, we have 
\(\Gamma\models\phi\,\Longrightarrow\,\Gamma\vdash\phi.\)
\end{theorem}
\begin{proof}
By the Compactness Theorem (\Cref{compactness}) we may assume that $\Gamma$ is a finite set and $\psi=\bigwedge\Gamma$.
Suppose $\psi\models\phi$. By Lemma \ref{WNF_PD_sound} we have 
\[\psi\equiv \bigbor_{\Omega\in \Lambda}\psi^\ast_\Omega~~\text{ and }~~\phi\equiv \bigbor_{\Delta\in \Lambda'}\phi^\ast_\Delta,\]
where $\Lambda$ and $\Lambda'$ are the (nonempty) sets of all maximal realizing sequences of $\psi$ and  $\phi$, respectively.


For each $\Omega\in \Lambda$, we have 
\(\psi_\Omega^\ast\models \bigbor_{\Delta\in \Lambda'}\phi^\ast_\Delta\). Since $\Omega$ and $\Delta$ are maximal realizing sequences,   $\psi_\Omega^\ast(p_{i_1},\dots,p_{i_n})$ and $\phi^\ast_\Delta(p_{i_1},\dots,p_{i_n})$ are classical formulas. By the standard argument these two classical formulas can be transformed into formulas in the classical disjunctive normal form in the deduction system of \CPL and therefore  by  \Cref{comp_dep_fr} also in the deduction system of  \PD, i.e., 
\begin{equation}\label{pd_cmpl_eq1}
\psi_\Omega^\ast\dashv\vdash\Theta_{X_\Omega}~~\text{ and }~~\phi^\ast_\Delta\dashv\vdash\Theta_{Y_\Delta},
\end{equation}
where $X_\Omega$, $Y_\Delta$ are some teams on $\{i_1,\dots,i_n\}$, and the formula $\Theta_X(p_{i_1},\dots,p_{i_n})$ is defined as in Lemma \ref{X_ThetaX}. By the Soundness Theorem we obtain that $\Theta_{X_\Omega}\models\bigbor_{\Delta\in \Lambda'}\Theta_{Y_\Delta}$. Now by the same argument as that in the proof of \Cref{PDbor_completeness} one derives that $\Theta_{X_\Omega}\vdash\Theta_{Y_{\Delta_0}}$ for some $\Delta_0\in \Lambda'$, which  implies that $\psi_\Omega^\ast\vdash\phi^\ast_{\Delta_0}$ by (\ref{pd_cmpl_eq1}). Hence we obtain by \Cref{si_derive} that $\psi^\ast_\Omega\vdash\phi$ for all $\Omega\in \Lambda$, which by \Cref{sek_derive} yields $\psi\vdash\phi$ as desired.
\end{proof}


\section{Concluding remarks} 


We have shown that there is a robust expressively complete propositional dependence logic for reasoning about dependencies on the propositional level. This logic can be fairly called {\em robust} because it has five different expressively equivalent formulations: \PD, \PDbor, \PID, \Inql and \PT, and each of these different formulations has its complete axiomatization. This logic adopts the team semantics, which extends the single-valuation semantics of classical propositional logic \CPL in a conservative manner: the team semantics reduces to the single-valuation semantics for formulas in the language of \CPL, and over singleton teams all connectives of the logic behave classically. Yet the move to \emph{teams} allows us to meaningfully characterize the dependencies between propositions, which cannot be expressed in \CPL. Also, two non-classical connectives (the intuitionistic disjunction and implication) are identified in the team semantics setting. Exactly with these different connectives (and atoms) we formed our different formulations of the logic.

These five formulations of the same logic are well-behaved also in the sense that all of them have disjunctive or conjunctive normal forms that resemble the same normal forms in \CPL. Making essential use of the disjunctive normal form of the logics, \cite{Yang15interpolation} proved the Craig's Interpolation Theorem for all of the five logics, and  \cite{IemhoffYang15} showed that the logics are  \emph{structurally complete} with respect to flat substitutions (i.e., the admissible rules  of these logics with respect to flat substitutions are derivable).


%
%
%
%
%
%

The study of dependence in logic originates from the first-order context. When quantifiers are added to our logics, many nice properties obtained in this paper fail immediately. As  already noted, first-order dependence logic is not axiomatizable, it is much weaker in expressive power than first-order intuitionistic dependence logic, and the equivalence given in \Cref{dep_v_df} does not apply to first-order dependence atoms. The results in this paper show that  the dependence notion on the propositional level is ``simple". This is justified by \Cref{expressive_pw} and the simplicity of the sound and complete deduction systems we defined in \Cref{sec:axiom}.

We finish by mentioning some relevant related areas. One area is the study of the computational complexity of propositional logics of dependence. A downward closed formula $\phi$ being satisfiable by some team is equivalent to $\phi$ being satisfiable by some singleton team. Since over singleton teams dependence atoms are always true and the logics have the same semantics as classical propositional logic \CPL, the satisfiability problem (SAT) for \PT and its fragments has the same complexity as SAT for \CPL, which is well-known to be NP-complete. As for the model checking problem (MC), \cite{eblo11} shows that MC for \PD and \PDbor are NP-complete, and it follows from \cite{MID_mc} that MC for \PID and \Inql are coNP-complete.  \PD is closely related to \emph{dependency quantified Boolean formulas}, whose validity problem is shown in \cite{Peterson2001957} to be NEXPTIME-complete. On the basis of this, \cite{Virtema2014} showed that the computational complexity of the validity problem of \PD is  NEXPTIME-complete.



Adding modalities to propositional dependence logic, one obtains {\em modal dependence logic}, which was introduced in \cite{VaMDL08}. This field has been very active in recent years. Another recent development is the rise of \emph{independence logic} (introduced in \cite{D_Ind_GV}), which is, {\em a priori}, stronger than dependence logic. The authors have a  paper  on \emph{propositional independence logic} under preparation. 



\section*{Acknowledgements}
The authors would like to thank Ivano Ciardelli, Pietro Galliani, Lauri Hella, Taneli Huuskonen, Rosalie Iemhoff,  Dick de Jongh, Juha Kontinen, Tadeusz Litak, Floris Roelofsen and Dag Westerst\aa hl for inspiring conversations related to this paper. We also thank the anonymous referee for reading the paper carefully and making many useful suggestions, which led to a great improvement in the presentation of our results.






\end{document}

\endinput